\documentclass[12pt,oneside,english]{amsart}

\usepackage{times}
\usepackage[T1]{fontenc}
\usepackage[latin1]{inputenc}
\usepackage{amssymb}

\makeatletter
 \theoremstyle{plain}    
 \newtheorem{thm}{Theorem}[section]
 \numberwithin{equation}{section} 
 \numberwithin{figure}{section} 
 \theoremstyle{plain}
 \theoremstyle{remark}
 \newtheorem*{rem*}{Remark}
 \theoremstyle{definition}
  \newtheorem{problem}[thm]{Problem}
 \theoremstyle{definition}
\newtheorem{prop-defn}[thm]{Proposition-Definition}
\theoremstyle{plain}
 \newtheorem{defn}[thm]{Definition}
 \theoremstyle{plain}    
 \newtheorem*{cor*}{Corollary}
 \theoremstyle{plain}    
 \newtheorem{conjecture}[thm]{Conjecture} 
 \theoremstyle{plain}    
 \newtheorem*{thm*}{Theorem} 
 \theoremstyle{plain}    
 \newtheorem{lem}[thm]{Lemma} 
 \theoremstyle{remark}    
 \newtheorem{notation}[thm]{Notation} 
 \theoremstyle{plain}    
 \newtheorem*{lem*}{Lemma} 
 \theoremstyle{plain}    
 \newtheorem{cor}[thm]{Corollary} 
 \theoremstyle{plain}    
 \newtheorem*{prop*}{Proposition} 
 \theoremstyle{remark}
 \newtheorem{rem}[thm]{Remark}
 \theoremstyle{plain}    
 \newtheorem{prop}[thm]{Proposition} 
 \usepackage{verbatim}
 \theoremstyle{definition}
 \newtheorem*{defn*}{Definition}

\def\Z{\mathbb {Z}}
\def\Q{\mathbb {Q}}
\def\R{\mathbb {R}}
\def\C{\mathbb {C}}
\def\bs{\backslash}

\newcommand{\wklim}[1]{\xrightarrow[#1]{\textrm{wk-*}}}

\newcommand{\ag}{\left(\mathfrak{a}\colon\!\mathfrak{g}\right)}
\newcommand{\hg}{\left(\mathfrak{h}_{\C}\colon\!\mathfrak{g}_{\C}\right)}

\newcommand{\rW}{W\!\ag}
\newcommand{\rroots}{\Delta\!\ag}
\newcommand{\cW}{W\!\hg}
\newcommand{\croots}{\Delta\!\hg}

\newcommand{\mk}{M\backslash K}
\newcommand{\V}{{L^2(\mk)}}

\newcommand{\Basis}{\mathcal{B}}

\newcommand{\Cic}{C^{\infty}_{\mathrm{c}}}
\newcommand{\D}{\Cic(X)}
\newcommand{\DK}{\D_K}
\newcommand{\DCi}{\Cic(C_i)}
\newcommand{\CX}{C^{\infty}(X)}
\newcommand{\Man}{Y}

\newcommand{\lieg}{\mathfrak{g}}
\newcommand{\lieh}{\mathfrak{h}}
\newcommand{\liea}{\mathfrak{a}}

\def\eqdef{\overset{\text{def}}{=}}
\def\isom{\simeq}
\def\restrict{\!\upharpoonright}
\def\subgp{<}

\DeclareMathOperator{\ad}{ad}
\DeclareMathOperator{\Ad}{Ad}
\DeclareMathOperator{\Ind}{Ind}
\DeclareMathOperator{\Lie}{Lie}
\DeclareMathOperator{\Op}{Op}
\DeclareMathOperator{\pr}{pr}
\DeclareMathOperator{\vol}{vol}

\def\PGL{\mathrm{PGL}}
\def\GL{\mathrm{GL}}
\def\PSL{\mathrm{PSL}}
\def\SL{\mathrm{SL}}
\def\SO{\mathrm{SO}}

\def\ie {i.e.\ }
\def\eg {e.g.\ }
\def\wrt{w.r.t.\ }

\usepackage{babel}

\makeatother
\begin{document}

\title{On Quantum unique ergodicity for locally symmetric spaces I}
\date{July 23rd, 2004}

\author[L. Silberman]{Lior Silberman}
\address{Lior Silberman\\
  Department of Mathematics\\
  Princeton University\\
  Princeton, NJ 08544-0001\\
  USA.
}
\email{lior@Math.Princeton.EDU}

\author[A. Venkatesh]{Akshay Venkatesh}
\address{Akshay Venkatesh\\
  Department of Mathematics\\
  Massachusetts Institute of Technology\\
  Cambridge, MA 02139-4307\\
  USA.
}
\thanks{The second author was supported in part by NSF Grant DMS-0245606.
Part of this work was performed at the Clay Institute
Mathematics Summer School in Toronto.}
\email{akshayv@Math.MIT.EDU}


\begin{abstract}
We construct an equivariant microlocal lift for locally symmetric spaces. 
In other words, we demonstrate how to lift, in a ``semi-canonical''
fashion, limits of eigenfunction measures on locally symmetric spaces to
Cartan-invariant measures on an appropriate bundle. The construction uses 
elementary features of the representation theory of semisimple real Lie
groups, and can be considered a generalization of Zelditch's results
from the upper half-plane to all locally symmetric spaces of noncompact type. 
This will be applied in a sequel to settle a version of the quantum unique
ergodicity problem on certain locally symmetric spaces.
\end{abstract}

\maketitle

\section{Introduction}

\subsection{General starting point: the semi-classical limit on Riemannian
manifolds}

Let $\Man$ be a compact Riemannian manifold, with the associated Laplace
operator $\Delta$ and Riemannian measure $d\rho$. An important problem
of harmonic analysis (or mathematical physics) on $\Man$ is understanding
the behaviour of eigenfunctions of $\Delta$ in the large
eigenvalue limit. The equidistribution problem asks whether for an
eigenfunction $\psi$ with a large eigenvalue $\lambda$, $\left|\psi(x)\right|$
is approximately constant on $\Man$. 
This can be approached ``pointwise'' and ``on average'' (bounding
$\left\Vert \psi\right\Vert _{L^{\infty}(\Man)}$ and
$\left\Vert \psi\right\Vert _{L^{p}(\Man)}$
in terms of $\lambda$, respectively), or ``weakly'': asking whether
as $\left|\lambda\right|\to\infty$, the probability measures defined
by $d\bar{\mu}_{\psi}(x)=\left|\psi(x)\right|^{2}d\rho(x)$ converge
in the weak-{*} sense to the ``uniform'' measure $\frac{d\rho}{\vol(\Man)}$.
For example, Sogge \cite{Sogge:LpBds} derives $L^{p}$ bounds for
$2\leq p\leq\infty$, and in the special case of Hecke eigenfunctions on
hyperbolic surfaces, Iwaniec and Sarnak \cite{IwaniecSarnak:LinftyBds}
gave a non-trivial $L^{\infty}$ bound. Here we will consider the weak-{*}
equidistribution problem for a special class of manifolds and eigenfunctions.

A general approach to the weak-{*} equidistribution problem was found
by \v Snirel$'$man \cite{Schnirelman:Avg_QUE}. To an eigenfunction
$\psi$ he associates a distribution $\mu_{\psi}$ on the unit cotangent
bundle $S^{*}\Man$ projecting to $\bar{\mu}_{\psi}$ on $\Man$.
This construction (the ``microlocal lift'') proceeds using the
theory of pseudo-differential operators and has the property
that, for any sequence
$\left\{\psi_{n}\right\} _{n=1}^{\infty}\subset L^{2}(\Man)$ with eigenvalues
$\lambda_{n}$ tending to infinity, any weak-{*} limit of the
$\mu_{n}=\mu_{\psi_{n}}$ is a probability measure on the unit tangent bundle
$S^{*}\Man$, invariant under the geodesic flow. Since any weak-{*} limit of
the $\mu_{n}$ projects to a weak-{*} limit of the $\bar{\mu}_{n}$, it suffices
to understand these limits; Liouville's measure $d\lambda$ on $S^{*}\Man$
plays here the role of the Riemannian measure on $\Man$.

This construction has a natural interpretation from the point of view
of semi-classical physics.  The geodesic flow on $\Man$ describes the motion
of a free particle (``billiard ball'').  $S^{*}\Man$ is (essentially) the
\emph{phase space} of this system, \ie the state space of the motion.
In this setting one calls a function $g\in C^\infty(S^{*}\Man)$ an
\emph{observable}.
The state space of the quantum-mechanical billiard is $L^2(\Man)$, with
the infinitesimal generator of time evolution $-\Delta$.  ``Observables'' here
are bounded self-adjoint operators $B\colon L^2(\Man)\to L^2(\Man)$. Decomposing
a state $\psi\in L^2(\Man)$ \wrt the spectral measure of $B$ gives a
probability measure on the spectrum of $B$ (which is the set of possible
``outcomes'' of the measurement).  The expectation value of the 
``measuring $B$ while the system is in the state $\psi$'' is then given
by the matrix element $\langle B\psi,\psi\rangle$.  In the particular case where
$B$ is a pseudo-differential operator with symbol $g\in C^\infty(S^{*}\Man)$, we
think of $B$ as a ``quantization'' of $g$, and any such a $B$ will be denoted
$\Op(g)$.  

We can now describe \v Snirel$'$man construction: it is given by
$\mu_\psi(g) = \langle \mathrm{Op}(g)\psi,\psi\rangle$.
This indeed lifts
$\bar{\mu}_\psi$, since for $g\in C^\infty(\Man)$ we can take $\Op(g)$ to be
multiplication by $g$.  If the $\psi$ are taken to be eigenfunctions then,
asymptotically, this construction does not depend on the choice
of ``quantization scheme,'' that is to say,
on the choice of the assignment $g\mapsto\Op(g)$.
Indeed, if $B_1, B_2$ have
the same symbol of order 0, and $-\Delta\psi=\lambda\psi$ 
(\ie ``$\psi$ is
an eigenstate of energy $\lambda$'') 
then one has
$\langle(B_1-B_2)\psi,\psi\rangle = O(\lambda^{-1/2})$. 

On a philosophical level we expect that at the limit of large energies, our
quantum-mechanical description to approach the classical one. We will not
formalize this idea (the ``correspondence principle''), but depend on
it for motivating our main question, whether ergodic properties of the
classical system persist in the semi-classical limit of the ``quantized''
version:

\begin{problem} (Quantum Ergodicity)
Let $\left\{ \psi_{n}\right\} _{n=1}^{\infty}\subset L^{2}(\Man)$
be an orthonormal basis consisting of eigenfunctions of the Laplacian.
\begin{enumerate}
\item What measures occur as weak-{*} limits of the $\left\{ \bar{\mu}_{n}\right\} $?
In particular, when does $\bar{\mu}_{n}\wklim{n\to\infty}d\rho$ hold?
\item What measures occur as weak-{*} limits of the $\left\{ \mu_{n}\right\} $?
In particular, when does $\mu_{n}\wklim{n\to\infty}d\lambda$ hold?
\end{enumerate}
\end{problem}
\begin{defn}
Call a measure $\mu$ on $S^{*}\Man$ a (microlocal) \emph{quantum limit}
if it is a weak-{*} limit of a sequence of distributions $\mu_{\psi_{n}}$
associated, via the microlocal lift, to a sequence of eigenfunctions
$\psi_{n}$ with $\left|\lambda_{n}\right|\to\infty$.
\end{defn}
In this language, the main problem is classifying the quantum limits
of the classical system, perhaps showing that the Liouville measure
is the unique quantum limit. As formalized by Zelditch
\cite{Zelditch:SL2_Lift_II} (for surfaces of constant negative curvature)
and Colin de Verdi{\`e}re \cite{CdV:Avg_QUE} (for general $\Man$),
the best general result known is still:

\begin{thm}
\label{thm: QUE_on_avg}Let $\Man$ be a compact manifold, $\left\{
\psi_{n}\right\} _{n=1}^{\infty}\subset L^{2}(\Man)$
an orthonormal basis of eigenfunctions of $\Delta$, ordered
by increasing eigenvalue. Then:
\end{thm}
\begin{enumerate}
\item (Weyl's law; see \eg \cite{Hormander:WeylLaw}) $\frac{1}{N}\sum_{n=1}^{N}\mu_{n}\wklim{N\to\infty}d\lambda$
holds with no further assumptions.
\item (\v Snirel$'$man-Zelditch-Colin de Verdi{\`e}re) Under the additional
assumption that the geodesic flow on $S^{*}\Man$ is ergodic, there exists
a subsequence $\left\{ n_{k}\right\} _{k=1}^{\infty}$ of density
$1$ s.t. $\mu_{n_{k}}\wklim{k\to\infty}d\lambda$.
\end{enumerate}
\begin{cor*}
For this subsequence, $\bar{\mu}_{n_{k}}\wklim{k\to\infty}d\rho$.
\end{cor*}

It was proved by Hopf \cite{Hopf:NegCurvErg} that the geodesic flow on a
manifold of negative sectional curvature is ergodic.  In this case, Rudnick
and Sarnak \cite{RudnickSarnak:Conj_QUE} conjecture a simple situation:

\begin{conjecture}
\label{con: RS_QUE}(Quantum unique ergodicity) Let $\Man$ be a compact
manifold of strictly negative sectional curvature. Then:
\begin{enumerate}
\item (QUE on $\Man$) $\bar{\mu}_{n}$ converge weak-{*} to the Riemannian
measure on $\Man$.
\item (QUE on $S^{*}\Man$) $d\lambda$ is the unique quantum limit on $\Man$.
\end{enumerate}
\end{conjecture}

We remark that \cite{RudnickSarnak:Conj_QUE} also gives
an example of a hyperbolic $3$-manifold $Y$,
a point $P \in Y$, and a sequence of eigenfunctions $\psi_n$
with eigenvalues $\lambda_n$,
such that $|\psi_n(P)| \gg \lambda_n^{1/4 - \epsilon}$.
The point $P$ is a fixed point of many Hecke
operators, and behaves 
in a similar fashion to the poles of a surface of revolution.
This remarkable phenomenon does not seem to contradict
Conjecture \ref{con: RS_QUE}. In the sequel to this paper
the scarcity of such points and their higher-dimensional analogues
will play an important role. 

\subsection{Past work: Quantum
unique ergodicity on hyperbolic surfaces and $3$-manifolds}

The quantum unique ergodicity question for hyperbolic surfaces
has been intensely investigated
over the last two decades. We recall some important results. 

Zelditch's work \cite{Zelditch:SL2_Lift_I,Zelditch:SL2_Lift_A_inv}
on the case of compact surfaces $\Man$ of constant negative curvature provided
a representation-theoretic alternative to the original construction of
the microlocal lift via the theory of pseudo-differential operators.
It is well-known that the universal cover of such a surface $\Man$ is
the upper half-plane $\mathbb{H}\isom \PSL_{2}(\R)/\SO_{2}(\R)$, so
$\Man=\Gamma\backslash\mathbb{H}$ for a uniform lattice
$\Gamma< G = \PSL_{2}(\R)$.  Then the $\SO_{2}(\R)\isom S^{1}$ bundle
$X=\Gamma\backslash \PSL_{2}(\R)\twoheadrightarrow \Man$
is isomorphic to the unit cotangent bundle of $\Man$. In this parametrization,
the geodesic flow on $S^{*}\Man$ is given by the action of the maximal
split torus
$A=\left\{ \left(\begin{array}{cc}
e^{t/2} \\
 & e^{-t/2}
\end{array}\right)\right\}$
on $X$ from the right. Zelditch's explicit microlocal lift starts
with the observation that an eigenfunction $\psi_{n}$ (considered
as a $K$-invariant function on $X$) can be thought of as the spherical
vector $\varphi_{0}^{(n)}$ in an irreducible $G$-subrepresentation
of $L^{2}(X)$. He then constructs another ({}``generalized'') vector
in this subrepresentation, a distribution $\delta^{(n)}$, 
and shows that the distribution given by $\mu_{\psi_{n}}(g)= \delta^{(n)}
(g\overline{\varphi_{0}^{(n)}})$
for $g\in \Cic(X)$ agrees (up to terms which decay as the $\lambda_{n}$
grow) with the microlocal lift. He then observes that the distribution
$\mu_{\psi_n}$ is exactly annihilated
by a differential operator of the form $H+\frac{J}{r_{n}}$ where
$H$ is the infinitesimal generator of the geodesic flow, $J$ a certain
(fixed) second-order differential operator, and $\lambda_{n}=-\frac{1}{4}-r_{n}^{2}$.
It is then clear that any weak-{*} limit taken as
$\left|\lambda_{n}\right|\to\infty$ will be annihilated (in the sense of
distributions) by the differential operator $H$, or in other words be
invariant under the geodesic flow.
Wolpert \cite{Wolpert:SL2_Lift_Fejer} made Zelditch's approach self-contained
by showing that the limits are positive measures without using
pseudo-differential calculus. For a clear exposition of the Zelditch-Wolpert
microlocal lift see \cite{Lindenstrauss:HH_QUE}.

Lindenstrauss's paper \cite{Lindenstrauss:HH_QUE} considers the case
of $\Man=\Gamma\bs\left(\mathbb{H}\times\cdots\times\mathbb{H}\right)$
for an irreducible lattice $\Gamma$ in
$\PSL_{2}(\R)\times\cdots\times \PSL_{2}(\R)$.
The natural candidates for $\psi_{n}$ here are not eigenfunctions
of the Laplacian alone, but rather of all the ``partial'' Laplacians
associated to each factor separately. 
Set now $G=\PSL_{2}(\R)^{h}$, $K=\SO_{2}(\R)^{h}$, $X=\Gamma\bs G$, $Y=\Gamma\bs G/K$,
and take $\Delta_{i}$ to be the Laplacian operator associated with the $i$th factor
(so that $\C\left[\Delta_{1},\ldots,\Delta_{h}\right]$ is the ring
of $K$-bi-invariant differential operators on $G$).
Assume that $\Delta_{i}\psi_{n}+\lambda_{n,i}\psi_{n}=0$,
where $\lim_{n\to\infty}\lambda_{n,i}=\infty$ for each $1\leq i\leq h$
separately. Generalizing the Zelditch-Wolpert construction, Lindenstrauss
obtains distributions $\delta^{(n)}\overline{\varphi_{0}^{(n)}}$
on $X$, projecting to $\bar{\mu}_{\psi_{n}}$on $Y$, and so that every
weak-{*} limit of these (a {}``quantum limit'') is a finite positive
measure invariant under the action of the full maximal split torus
$A^{h}$. He then proposes the following version of QUE, also due
to Sarnak:

\begin{problem}
\label{pro: SQUE-Y}(QUE on locally symmetric spaces) Let $G$ be
a connected semi-simple Lie group with finite center. Let $K<G$ be
a maximal compact subgroup, $\Gamma<G$ a lattice, $X=\Gamma\bs G$,
$Y=\Gamma\bs G/K$.
Let $\left\{ \psi_{n}\right\} _{n=1}^{\infty}\subset L^{2}(Y)$
be a sequence of normalized eigenfunctions of the ring of 
$G$-invariant differential operators on $G/K$, with the eigenvalues \wrt
the Casimir operator tending to $\infty$ in absolute value.
Is it true that $\bar{\mu}_{\psi_{n}}$ converge
weak-{*} to the normalized projection of the Haar measure to $Y$?
\end{problem}

\subsection{This paper: Quantum unique ergodicity on locally symmetric
spaces}

This paper is the first of two papers on this general problem. The
main result of the present paper (Theorem \ref{mainthm} below) is the
construction of the microlocal lift in this setting. 
We will impose a mild non-degeneracy condition on the sequence of
eigenfunctions (see Section \ref{sub: positivity}; the assumption
essentially amounts to asking
that all eigenvalues tend to infinity, at the same rate for operators
of the same order.)

With $K$ and $G$ as in Problem \ref{pro: SQUE-Y},
let $A$ be as in the Iwasawa decomposition $G = NAK$, \ie 
$A = \mathrm{exp}(\mathfrak{a})$ where $\mathfrak{a}$ is a maximal
abelian subspace of $\mathfrak{p}$. (Full definitions are given
in Section \ref{sub: notation}). 
For $G = \SL_n(\R)$ and $K = \SO_n(\R)$,
one may take $A$ to be the
subgroup of diagonal matrices with positive entries. 
Let $\pi \colon X \to Y$ be the projection.
We denote by $dx$ the $G$-invariant probability measures 
on $X$, and by $dy$ the projection of this measure to $Y$. 

The content of the Theorem that follows amounts, roughly,
to a ``$G$-equivariant microlocal lift'' on $Y$.

\begin{thm} \label{mainthm}
Let $\left\{ \psi_{n}\right\} _{n=1}^{\infty}\subset L^{2}(Y)$
be a non-degenerate sequence of normalized eigenfunctions,
whose eigenvalues approach $\infty$. 
Then, after replacing $\psi_n$ by an appropriate subsequence, 
there exist functions $\tilde{\psi}_n \in L^2(X)$ 
and distributions $\mu_n$ on $X$ such that:
\begin{enumerate}
\item \label{thm:claim1} The projection of $\mu_n$ to $Y$
coincides with $\bar{\mu}_n$, \ie 
$\pi_{*} \mu_n = \bar{\mu}_n$. 

\item \label{thm:claim2} Let $\sigma_n$ be the measure $|\tilde{\psi}_n(x)|^2 dx$ 
on $X$. 
Then, for every $g \in \D$, we have
$\lim_{n \to \infty} (\sigma_n(g) - \mu_n(g)) = 0$. 

\item \label{thm:claim3} Every weak-{*} limit $\sigma_{\infty}$ of the measures
$\sigma_n$ (necessarily a positive measure of mass $\leq 1$)
is $A$-invariant. 
\item \label{prop:eq} (Equivariance). Let $E\subset\mathrm{End}_G(\CX)$ 
be a $\C$-subalgebra of bounded endomorphisms of $\CX$, commuting with the
$G$-action.
Noting that each $e \in E$ induces an endomorphism of $C^{\infty}(Y)$,
suppose that $\psi_n$ is an eigenfunction for $E$ (\ie 
$E  \psi_n \subset \C \psi_n$). 
Then we may
choose $\tilde{\psi}_n$
so that $\tilde{\psi}_n$ is an eigenfunction for $E$ with the same eigenvalues
as $\psi_n$, \ie for all $e \in E$ there exists $\lambda_e \in \C$ such that
$e \psi_n = \lambda_e \psi_n, e \tilde{\psi}_n = \lambda_e \tilde{\psi}_n$. 
	
\end{enumerate}
\end{thm}

We first remark that the distributions $\mu_n$ (resp.
the measures $\sigma_n$) generalize the constructions
of Zelditch (resp. Wolpert). Although, in view of
(\ref{thm:claim2}), they carry roughly equivalent information,
it is convenient to work with both simultaneously:
the distributions $\mu_n$ are canonically defined and easier
to manipulate algebraically, whereas the measures $\sigma_n$
are patently positive and are central to the arguments
in the sequel to this paper.

\begin{proof}
For simplicity, we first write the proof in detail for the case where
$G$ is simple (the modifications necessary in the general case are
discussed in Section \ref{sub: extend}).

In Section \ref{sub: def-lift}
we define the distributions $\mu_{n}$.
(In the language
of Definition \ref{def: lift}, we take 
$\mu_n = \mu_{\psi_n}(\varphi_0 , \delta)$). 

Claim (\ref{thm:claim1}) is established in
 Lemma \ref{lem:projects}. 

In Section \ref{sub: positivity} we introduce the non-degeneracy
condition. 
 Proposition \ref{prop: equivariance}
defines $\tilde{\psi}_n$ and
establishes the claims (\ref{thm:claim2}) and (\ref{prop:eq}). 
(Observe that this Proposition establishes (\ref{thm:claim2})
only for $K$-finite test functions $g$. Since the extension
to general $g$
is not necessary for any of our applications, we omit the proof.)

Finally, in section \ref{sub: A-inv} we
establish claim (\ref{thm:claim3}) (Corollary \ref{cor: A-inv})
by finding enough differential operators annihilating $\mu_{n}$.
\end{proof}

\begin{rem}\label{rem: mainthm}\ \\
\begin{enumerate}
\item
It is important to verify that non-degenerate sequences of
eigenfunctions exist.  In the co-compact case (\eg for the
purpose of Theorem \ref{thm:que}), it was shown in
\cite{DuistermaatKolkVaradarajan:WeylLaw,DuistermaatKolkVaradarajan:WeylLaw_Erratum}
that a positive proportion of the unramified spectrum lies in every
open subcone of the Weyl chamber (for definitions see
Theorem \ref{thm: ss-gps-theory} and the discussion in
Section \ref{subsec:motivation}). This is also
expected to hold for finite-volume \emph{arithmetic} quotients $Y$.
For example, \cite[Thm.\ 5.3]{SDMiller:SL3_WeylLaw} treats the
case of $\SL_3(\Z)\bs\SL_3(\R)/\SO_3(\R)$.

\item 
We shall use the phrase \emph{non-degenerate quantum limit} to denote
any weak-{*} limit of $\sigma_n$, where notations
are as in Theorem \ref{mainthm}.  Note that
if $\sigma_{\infty}$ is such a limit, then
claim (\ref{thm:claim2}) of the Theorem shows 
that there exists a subsequence $(n_k)$ of the integers
such that
$\sigma_{\infty}(g) =  \lim_{n_k \rightarrow \infty} \mu_{n_k}(g)$
for all $g \in \D$. 
Depending on the context, we shall therefore use the notation
$\sigma_{\infty}$ or $\mu_{\infty}$ for a non-degenerate quantum limit.

\item 
It is not necessary to pass to a subsequence in 
Theorem \ref{mainthm}. See Remark \ref{rem:subseq}. 

\item \label{rem: micro-or-not}
It is likely that the $A$-invariance
aspect of Theorem \ref{mainthm} could be established
by standard microlocal methods; however, 
the equivariance property does not follow readily from these methods
and is absolutely crucial in applications.
It will be used, in the sequel to this paper,
in the situation where $E$ is an algebra of endomorphisms
generated by Hecke correspondences. 

\item 
The measures $\sigma_{\infty}$ all are invariant
by the compact group $M = Z_K(\mathfrak{a})$. In
fact, Theorem \ref{mainthm} should strictly be interpreted
as lifting measures to $X/M$ rather than $X$.

\item
Theorem \ref{mainthm} admits a natural geometric interpretation. 
Informally, the bundle $X/M \rightarrow Y$ may be regarded as a bundle
parameterizing maximal flats in $Y$, and the $A$-action
on $X/M$ corresponds to ``translation along flats.'' 
We refer to Section \ref{chamber} for a further discussion
of this point. 

\end{enumerate}
\end{rem}

The existence of the microlocal lift already places a restriction
on the possible weak-{*} limits of the measures
$\left\{ \bar{\mu}_{n}\right\}$ on $Y$. 
In particular, Theorem \ref{mainthm} has the following corollary (in this
regard see also Remark \ref{rem: mainthm}(\ref{rem: micro-or-not})).
\begin{cor}  \label{cor:flat}
Let $\left\{ \psi_{n}\right\} _{n=1}^{\infty}\subset L^{2}(Y)$ be
a non-degenerate sequence of normalized eigenfunctions such that $\bar{\mu}_{\psi_{n}}$
converge in the weak-{*} topology to a limit measure $\bar{\mu}_{\infty}$.
Then $\bar{\mu}_{\infty}$ is the projection to $Y$ of an $A$-invariant
measure $\mu_{\infty}$ on $X$. In particular, the support of $\mu_{\infty}$
must be a union of maximal flats.
\end{cor}

More importantly, Theorem \ref{mainthm} allows us to pose a new version of
the problem:

\begin{problem}
\label{pro: SQUE-X}(QUE on homogeneous spaces) In the setting of
Problem \ref{pro: SQUE-Y}, is the $G$-invariant measure on $X$
the unique non-degenerate quantum limit?
\end{problem}

\subsection{Arithmetic QUE. Sequel to this paper.} 
The sequel to this paper will resolve Problem \ref{pro: SQUE-X}
for various higher rank symmetric spaces, 
in the context of \emph{arithmetic} quantum limits. We briefly
recall their definition and significance. 

Let $\Man$ be (for example) a negatively curved manifold. 
In general, we believe that the multiplicities of the Laplacian $\Delta$
acting on $L^2(\Man)$ are quite small, \ie 
the $\lambda$-eigenspace has dimension $\ll_{\epsilon} \lambda^{\epsilon}$. 
This question seems extremely difficult even for
$\SL_2(\Z) \bs \mathbb{H}$, and no better bound is known than the general
$O(\lambda^{1/2}/\log(\lambda))$, valid for all negatively curved manifolds.

However, even lacking information on the multiplicities,
it transpires that in many natural instances
we have
a \emph{distinguished basis} for $L^2(\Man)$. 
In that context, it is then natural
to ask whether Problem \ref{pro: SQUE-Y} or Problem 
\ref{pro: SQUE-X}
can be resolved with respect to this distinguished basis.
Since it is believed that the $\Delta$-multiplicities
are small, this modification is, philosophically, not too far from the
original question. However, it is in many natural cases far more tractable.

The situation of having (something close to)
a distinguished basis occurs for $\Man =\Gamma \backslash G/K$
and $\Gamma \subset G$ arithmetic. This distinguished
basis is obtaining by simultaneously diagonalizing
the action of Hecke operators. We shall not
give precise definitions here; in any case,
we refer to quantum limits arising
from subsequences of the distinguished basis
as \emph{arithmetic quantum limits}.

In the second paper we apply this results of this paper to the study
of arithmetic quantum limits. In particular we
settle the conjecture in the case where
$\Gamma$ arises from 
the multiplicative group of a division algebra of prime degree
over $\Q$. For brevity, we state the result in the language
of automorphic forms; in particular, $\mathbb{A}$ is the ring
of adeles of $\mathbb{Q}$.

\begin{thm} \label{thm:que}
(QUE for division algebras of prime degree) Let $\mathbb{D}/\Q$ be
a division algebra of prime degree $d$
and let $\mathbb{G}=\mathrm{P}\mathbb{D}^{\times}$
be the associated projective general linear group. Assume that $\mathbb{D}$
is unramified at $\infty$, \ie that
$\mathbb{G}(\R)\isom \PGL_{d}(\R)$.
Let $K_{f}<\mathbb{G}(\mathbb{A}_{f})$ be an open compact subgroup,
and let $\Gamma<\mathbb{G}(\R)$ be the (congruence) lattice such
that $X=\Gamma\backslash\mathbb{G}(\R)\isom\mathbb{G}(\Q)\backslash\mathbb{G}(\mathbb{A})/K_{f}$.
Then the normalized Haar measure is the unique non-degenerate arithmetic
quantum limit on $X$.
\end{thm}

We expect the techniques developed for the proof of Theorem \ref{thm:que}
will generalize at least to some other locally symmetric spaces,
the case of $\mathbb{D}$ being the simplest; but there are considerable
obstacles to obtaining a theorem for \emph{any} arithmetic locally symmetric
space at present. 

Let us make some remarks about the proof of Theorem \ref{thm:que}.
Our approach follows that of Lindenstrauss in \cite{Lindenstrauss:SL2_QUE}
which established the above theorem for division algebras of degree $2$.
This approach is based on result toward the classification of the
$A$-invariant measures on $X$.  To apply such a result one needs to show
further regularity of the limit measure -- that $A$ acts on every
$A$-ergodic component of $\mu_{\infty}$ with positive entropy.
This was proved for $G = \mathrm{SL}_2$ by Bourgain and Lindenstrauss 
in \cite{LindenstraussBourgain:SL2_Ent}.
In the higher-rank case we rely on recent results toward the classification
of the $A$-invariant measures on $X$,
due to Einsiedler-Katok \cite{EinsiedlerKatok:SLn_Rigid}, and prove
the positive entropy property of $\mu_\infty$.

Establishing positive entropy in higher rank is quite involved.
The equivariance (property (\ref{prop:eq})
of Theorem \ref{mainthm}), applied with $E$ the Hecke algebra,
plays a crucial role, just as in \cite{LindenstraussBourgain:SL2_Ent}. 
The proof utilizes
a study of the behavior of eigenfunctions on Bruhat-Tits buildings
and consideration of certain Diophantine questions (these questions
are higher-rank versions of the questions: to what extent can CM
points of bounded height cluster together?)

\subsection{Acknowledgments}
We are deeply indebted both to P. Sarnak and to E. Lindenstrauss for
their encouragement. We have also benefited from discussions
with J. Bernstein, A. Reznikov and E. Lapid. 

\section{Notation} \label{sec: Notation}

Section \ref{sub: notation} defines mostly standard notation and terminology
pertaining to semisimple groups and their root systems
(we generally follow \cite{Knapp:Rpn_Th_SS_Gps}). 
Section \ref{sub: repn-theory} sets up the basic theory of spherical
representations; the reader
may wish to read at least Definitions \ref{def: vdef} and \ref{defn:repdefn}.
Section \ref{sub: functional} defines the various function spaces
we will have need of; the notation here is fairly standard. 

\subsection{\label{sub: notation}General notation}

Let $G$ denote a non-compact connected simple Lie group with finite
center (we discuss generalizations to this in Section \ref{sub: extend}).
We choose a Cartan involution $\Theta$ for $G$, and let $K<G$ be
the $\Theta$-fixed maximal compact subgroup. Let $S=G/K$ be the
symmetric space, with $x_{K}\in S$ the point with stabilizer $K$.
We fix a $G$-invariant metric on $S$. 
To normalize it, we observe that the tangent
space at the point $x_{K}\in S$ is identified with $\mathfrak{p}$
(see below), and we endow it with the Killing form.

For a lattice $\Gamma<G$ we set $X=\Gamma\bs G$ and $Y=\Gamma\bs G/K$,
the latter being a locally symmetric space of non-positive curvature.
We normalize the Haar measures $dx$ on $X$, $dk$ on $K$ and $dy$ on $Y$
to have total mass $1$ (here $dy$ is the pushforward of $dx$ under the
the projection from $X$ to $Y$ given by averaging \wrt $dk$).

Let $\mathfrak{g}=\Lie(G)$, and let $\theta$ denote the
differential of $\Theta$, giving the Cartan decomposition $\mathfrak{g}=\mathfrak{k}\oplus\mathfrak{p}$
with $\mathfrak{k}=\Lie(K)$. Fix now a maximal abelian subalgebra
$\mathfrak{a}\subset\mathfrak{p}$. 

We denote by  $\mathfrak{a}_{\C}$ the complexification $\mathfrak{a}
\otimes_{\R} \C$;  we shall occasionally write $\mathfrak{a}_{\R}$
for $\mathfrak{a}$ for emphasis in some contexts.
We denote by $\mathfrak{a}^{*}$ (resp. $\mathfrak{a}_{\C}^{*}$) 
the real dual (resp. the complex dual) of $\mathfrak{a}$;
again, we shall occasionally write $\mathfrak{a}^{*}_{\R}$ for
$\mathfrak{a}^{*}$.  For $\nu \in \mathfrak{a}_{\C}^{*}$,
we define $\mathrm{Re}(\nu), \mathrm{Im}(\nu) \in \mathfrak{a}_{\R}^{*}$
to be the real and imaginary parts of $\nu$, respectively. 

For $\alpha\in\mathfrak{a}^{*}$
set $\mathfrak{g}_{\alpha}=\left\{ X\in\mathfrak{g}\right|\forall H\in\mathfrak{a}:\ad(H)X=\alpha(H)X\left.\right\} $,
$\rroots=\left\{ \alpha\in\mathfrak{a}^{*}\setminus\left\{ 0\right\} \mid\mathfrak{g}_{\alpha}\neq\left\{ 0\right\} \right\} $
and call the latter the (restricted) \emph{roots} of
$\mathfrak{g}$ \wrt $\mathfrak{a}$.
The subalgebra $\mathfrak{g}_{0}$ is $\theta$ invariant, and hence
$\mathfrak{g}_{0}=(\mathfrak{g}_{0}\cap\mathfrak{p})\oplus(\mathfrak{g}_{0}\cap\mathfrak{k})$.
By the maximality of $\mathfrak{a}$ in $\mathfrak{p}$, we must then
have $\mathfrak{g}_{0}=\mathfrak{a}\oplus\mathfrak{m}$
where $\mathfrak{m}=Z_{\mathfrak{k}}(\mathfrak{a})$.

The Killing form of $\mathfrak{g}$ induces a standard inner product
$\left\langle \cdot,\cdot\right\rangle $ on $\mathfrak{a}^{*}$ \wrt 
which $\rroots\subset\mathfrak{a}^{*}$ is a root system. The associated
Weyl group, generated by the root reflections $s_{\alpha}$, will
be denoted $\rW$. This group is also canonically isomorphic to the
analytic Weyl groups $N_{G}(A)/Z_{G}(A)$ and $N_{K}(A)/Z_{K}(A)$.
The fixed-point set of any $s_{\alpha}$ is a hyperplane in $\mathfrak{a}^{*}$,
called a \emph{wall}.
The connected components of the complement of
the union of the walls are cones, called the \emph{(open) Weyl chambers}.
A subset $\Pi \subset \rroots$ will be called a \emph{system of simple roots}
(by abuse of notation a ``simple system'') if every root can be uniquely 
expressed as an integral combination of elements of $\Pi$
with either all coefficients non-negative
or all coefficients non-positive.
For a simple system $\Pi$, the open cone $C_{\Pi}=\left\{ \nu\in\mathfrak{a}^{*}\mid\forall\alpha\in\Pi:\left\langle \nu,\alpha\right\rangle >0\right\} $
is an open Weyl chamber, and the map $\Pi\mapsto C_{\Pi}$ is a $1-1$
correspondence between simple systems and chambers. The Weyl group
acts simply transitively on the chambers and simple systems. The closure
of an open chamber will be called a closed chamber. The action of
$\rW$ on $\mathfrak{a}^{*}$ extends in the complex-linear way to
an action on $\mathfrak{a}_{\C}^{*}$
preserving $i\mathfrak{a}^{*}\subset\mathfrak{a}_{\C}^{*}$, and we
call an element $\nu\in \mathfrak{a}_{\C}^{*}$ \emph{regular} if it
is fixed by no $w\in\rW$. We use $\rho=\frac{1}{2}\sum_{\alpha>0}(\dim\mathfrak{g}_{\alpha})\alpha\in\mathfrak{a}^{*}$
to denote half the sum of the positive (restricted) roots.

Fixing a simple system $\Pi$ we get a notion of positivity. For $\mathfrak{n}=\oplus_{\alpha>0}\mathfrak{g}_{\alpha}$
and $\bar{\mathfrak{n}}=\Theta\mathfrak{n}$ we have $\mathfrak{g}=\mathfrak{n}\oplus\mathfrak{a}\oplus\mathfrak{m}\oplus\bar{\mathfrak{n}}$
and (Iwasawa decomposition)
$\mathfrak{g}=\mathfrak{n}\oplus\mathfrak{a}\oplus\mathfrak{k}$.
By means of the Iwasawa decomposition, we may uniquely write every $X\in\mathfrak{g}$ in the form $X=X_{\mathfrak{n}}+X_{\mathfrak{a}}+X_{\mathfrak{k}}$.
We sometimes also write $H_{0}(X)$ for $X_{\mathfrak{a}}$.

Let $N,\, A<G$ be the subgroups corresponding to the subalgebras
$\mathfrak{n},\,\mathfrak{a}\subset\mathfrak{g}$ respectively, and
let $M=Z_{K}(\mathfrak{a})$. Then $A$ is a maximal split torus in
$G$, and $\mathfrak{m}=\Lie(M)$, though $M$ is not necessarily
connected. Moreover $P_{0}=NAM$ is a minimal parabolic subgroup of
$G$, with the map $N\times A\times M\to P_{0}$ being a diffeomorphism. The
map $N\times A\times K\to G$ is a (surjective) diffeomorphism (Iwasawa decomposition),
so for $g\in G$ there exists a unique $H_{0}(g)\in\mathfrak{a}$
such that $g=n\exp(H_{0}(g))k$ for some $n\in N$, $k\in K$. The
map $H_{0}:G\to\mathfrak{a}$ is continuous; restricted to $A$ it is the inverse
of the exponential map.

Let $\mathfrak{g}_{\C}=\mathfrak{g}\otimes_{\R}\C$ denote the complexification
of $\mathfrak{g}$. It is a complex semi-simple Lie algebra. Let $\theta_{\C}$ denote the \emph{complex-linear}
extension of $\theta$ to $\mathfrak{g}_{\C}$. It is \emph{not} a
\emph{}Cartan involution of $\mathfrak{g}_{\C}$. We fix a maximal
abelian subalgebra $\mathfrak{b}\subset\mathfrak{m}$ and set $\mathfrak{h}=\mathfrak{a}\oplus\mathfrak{b}$.
Then $\mathfrak{h}_{\C}=\mathfrak{h}\otimes\C\subset\mathfrak{g}_{\C}$
is a Cartan subalgebra, with the associated root system $\croots$
satisfying $\rroots=\left\{ \alpha\restrict_{\mathfrak{a}}\right\} _{\alpha\in\croots}\setminus\left\{ 0\right\} $.
Moreover, we can find a system of simple roots $\Pi_{\C}\subset\croots$
and a system of simple roots $\Pi \subset \rroots$
such that the positive roots \wrt $\Pi$ are precisely
the nonzero restrictions of the positive roots \wrt $\Pi_{\C}$. 
We fix such
a compatible pair of simple systems, and let $\rho_{\mathfrak{h}}$
denote half the sum of the roots in $\croots$, positive \wrt $\Pi_{\C}$.

Let $F_{0}\subset\croots$ consist of the roots
that restrict to $0$ on $\mathfrak{a}$, $F_{0}^{+}\subset F_{0}$
those positive \wrt $\Pi_{\C}$. Let $\mathfrak{n}_{M}=\oplus_{\alpha\in F_{0}^{+}}(\mathfrak{g}_{\C})_{\alpha}$,
$\bar{\mathfrak{n}}_{M}=\oplus_{\alpha\in F_{0}^{+}}(\mathfrak{g}_{\C})_{-\alpha}$.
Then $\mathfrak{m}_{\C}=\mathfrak{n}_{M}\oplus\mathfrak{b}_{\C}\oplus\bar{\mathfrak{n}}_{M}$
and
$\mathfrak{g}_{\C}=\mathfrak{n}_{\C}\oplus\mathfrak{n}_{M}\oplus\mathfrak{h}_{\C}\oplus\bar{\mathfrak{n}}_{M}\oplus\bar{\mathfrak{n}}_{\C}.$

For $\nu\in\mathfrak{a}_{\C}^{*}$, set $\left\Vert \nu\right\Vert ^{2}=\left\langle \mathrm{Re}(\nu),\mathrm{Re}(\nu)\right\rangle +\left\langle \textrm{Im}(\nu),\textrm{Im}(\nu)\right\rangle $
(with the inner products taken in $\mathfrak{a}_{\R}^{*})$.

If $\mathfrak{l}_{\C}$ is a complex Lie algebra, then we denote by
$U(\mathfrak{l}_{\C})$ its universal enveloping algebra,
and by $\mathfrak{Z}(\mathfrak{l}_{\C})$
its center. In particular we set $\mathfrak{Z}=\mathfrak{Z}(\mathfrak{g}_{\C})$.

\subsection{\label{sub: repn-theory}Spherical Representations and the model
$(V_{K},I_{\nu})$.}

We recall some facts from the representation theory of compact and
semi-simple groups. At the end of this section we analyze a model
(the ``compact picture'') for the spherical dual of $G$.

\begin{thm}
\label{thm: cpt-gps-theory}\cite[Th.\ 1.12]{Knapp:Rpn_Th_SS_Gps}
Let $K$ be a compact topological group
and let $\hat{K}_\mathrm{fin}$ be the set of equivalence classes
of irreducible finite-dimensional unitary representations of $K$. 
\begin{enumerate}
\item (Peter-Weyl) Every $\rho\in\hat{K}_\mathrm{fin}$ occurs discretely
in $L^{2}(K)$ with multiplicity equal to its dimension $d(\rho)$.
Moreover, $L^2(K)$ is isomorphic to the Hilbert direct sum of its
isotypical components $\{L^2(K)_\rho\}_{\rho\in\hat{K}_\mathrm{fin}}$.

\item Let $\pi:K\to \GL(W)$ be a representation of $K$ on the locally convex
complete space $W$. Then $\oplus_{\rho \in \hat{K}} W_{\rho}$ 
is dense in $W$, where $W_{\rho}$ is the $\rho$-isotypical subspace.

\item Every irreducible representation of $K$ on a locally convex,
complete space is finite-dimensional and hence unitarizable.
In particular, $\hat{K}_\mathrm{fin}$ is the unitary dual of $K$.

\item For $K$ as in Section \ref{sub: notation}, $\hat{K}$ is countable.
\end{enumerate}
\end{thm}
Note that while \cite[Th.\ 1.12(c-e)]{Knapp:Rpn_Th_SS_Gps}
are only claimed for unitary representations on Hilbert spaces,
their proofs only rely on the action of the convolution algebra $C(K)$
on representations of $K$, and hence carry over with
little modification to the more general context needed here.
The last conclusion follows from the separability of $L^2(K)$,
which in turn follows from the separability of $K$.
\begin{notation}
Let $\pi:K\to \GL(W)$ be as above. The algebraic direct sum
$W_{K}\eqdef\oplus_{\rho\in\hat{K}}W_{\rho}$
consists precisely of these $w\in W$ which generate a finite-dimensional
$K$-subrepresentation. We refer to $W_K$ as the space of $K$-finite vectors. 
We will use $W^{K}$ to denote these vectors of $W$ fixed by $K$.
\end{notation}

\begin{defn} \label{def: vdef}Set $V = \V$, and set $V_K \subset V$
to be the space of $K$-finite vectors. Let $C^{\infty}(\mk)$ be the
smooth subspace, $C^{\infty}(\mk)'$ the space of distributions on $\mk$.
Let $V_K'$ (resp. $V'$) be the dual to $V_K$ (resp. $V$). 
Then we have natural inclusions
$V_K \subset C^{\infty}(\mk) \subset V$ and
$V_K' \supset C^{\infty}(\mk)'  \supset V'$; further, we have
(Riesz representation) a conjugate-linear isomorphism
\begin{equation} \label{eq:inclusion} V 
\stackrel{T}{\hookrightarrow} V' \end{equation}
where the map $T\colon V \to V'$ is defined via the rule
$T(f)(g) = \langle g,f \rangle_V = \int_{\mk} g \overline{f} dk$. 

Fix an increasing exhaustive
sequence of finite dimensional $K$-stable subspaces of $V_K$, \ie a sequence
$V_1 \subset V_2 \subset \dots \subset V_N \subset V_{N+1} \subset \dots$
of subspaces such that $\cup_{i=1}^{\infty} V_i = V_K$
and each $V_i$ is a $K$-subrepresentation.

For $\Phi \in V_K'$ and $1 \leq N \in \mathbb{Z}$, define the $N$-truncation
of $\Phi$ as the unique element $\Phi_N \in V_N$
such that $T(\Phi_N) - \Phi$ annihilates $V_N$. 

Finally let $\varphi_{0} \in V_K$ be the function that is identically $1$. 
\end{defn}

\begin{defn} \label{def:deltaseq}
Let $\mu$ be a regular Borel measure on a space $X$.  Call a sequence
of non-negative functions $\left\{ f_{j}\right\} \in L^{1}(\mu)$
a $\delta$-\emph{sequence} at $x \in X$ if, for every $j$,
$\int f_{j}d\mu=1$, and moreover if, for every $g\in C(X)$,
$\lim_{j\to\infty}\int f_{j}\cdot gd\mu=g(x)$.
\end{defn}
\begin{lem}
\label{lem: sqrt-delta}There exists a sequence
$\left\{ f_{j}\right\} _{j=1}^{\infty}\subset V_{K}$
such that $|f_{j}|^{2}$ is a $\delta$-sequence on $\mk$.
\end{lem}
\begin{proof}
Let $\left\{ h_{j}\right\} _{j=1}^{\infty}\subset C(\mk)$ be a
$\delta$-sequence.
By the Peter-Weyl theorem $V_{K}$ is dense in $C(\mk)$, so that
for every $j$ we can choose $f'_{j}\in V_{K}$ such that
$\left\Vert \sqrt{h_{j}(k)}-f_{j}'(k)\right\Vert _{\infty}\leq\frac{1}{2^{j}}$.
Then one may take $f_{j}=\frac{f'_{j}}{\left\Vert f_{j}'\right\Vert _{2}}$.
\end{proof}
Secondly, we recall the construction of the spherical principal series
representations of a semi-simple Lie group. An irreducible representation
of $G$ is \emph{spherical} if it contains a 
$K$-fixed vector. Such a vector is necessarily unique up to scaling. 

To any $\nu\in\mathfrak{a}_{\C}^{*}$ we associate the character $\chi_{\nu}(p)=\exp(\nu(H_{0}(p))$
of $P_{0}$ and the induced representation with $\left(\mathfrak{g},K\right)$-module
\begin{equation} \label{eq:indrep}\Ind_{P_{0}}^{G}\chi_{\nu}=\left\{ f\in{C^{\infty}(G)}_{K}\mid\forall p\in
P,g\in G:f(pg)=e^{\left\langle \nu+\rho,H_{0}(p)\right\rangle }f(g)\right\}
\end{equation}
By the Iwasawa decomposition, every $f\in\Ind_{P_{0}}^{G}\chi_{\nu}$
is determined by its restriction to $K$;
this restriction defines an element of the space
$V_{K}$. Conversely, every $f\in V_{K}$ extends uniquely to a member of
$\Ind_{P_{0}}^{G}\chi_{\nu}$.

\begin{defn} \label{defn:repdefn}
For $\nu \in \mathfrak{a}_{\C}^{*}$, we denote by
$(I_{\nu}, V_K)$ the representation of $\mathfrak{g}$
on $V_K$ fixed by the discussion above;
we shall also use $I_{\nu}$ to denote the
corresponding action of $\mathfrak{g}$ on $C^{\infty}(\mk)$
and of $G$ on $V$. 
We shall denote by $I_{\nu}'$ the dual action of $\mathfrak{g}$
on either $V_K'$ or $C^{\infty}(\mk)'$. 
\end{defn}

Note also that $\varphi_{0} \in V_K$ (see Definition \ref{def: vdef})
is a spherical vector for the representation $(I_{\nu}, V_K)$.

\begin{thm}
\label{thm: ss-gps-theory}(The unitary spherical dual;
references are drawn from \cite{Knapp:Rpn_Th_SS_Gps})
\begin{enumerate}

\item\
For any $\nu \in \liea_{\C}^{*}$, $\Ind_{P_{0}}^{G}\chi_{\nu}$
has a unique spherical irreducible subquotient, to be denoted $\pi_{\nu}$.
[Th.\ 8.37] Any spherical irreducible unitary
representation of $G$ is isomorphic to $\pi_{\nu}$ for
some $\nu$.  [Th.\ 8.38] We have $\pi_{\nu_{1}}\isom\pi_{\nu_{2}}$
iff there exists $w\in\rW$ such that $\nu_{2}=w\nu_{1}$.

\item\
[\S 7.1-3] If $\mathrm{Re}(\nu)=0$ then $\Ind_{P_{0}}^{G}\chi_{\nu}$
is unitarizable, with the invariant Hermitian form given
by $\left\langle f,g\right\rangle =\int_{\mk}f(k)\overline{g(k)}dk$.
This representation has a unique spherical summand
(necessarily isomorphic to $\pi_{\nu}$), and we let
$j_{\nu}\colon V_K\to\pi_{\nu}$ denote the orthogonal projection map.
[Th.\ 7.2] If $\nu$ is regular then $\Ind_{P_{0}}^{G}\chi_{\nu}$ is irreducible.

\item\
[\S 16.5(7) \& Th. 16.6] If $\pi_{\nu}$ is unitarizable then $\mathrm{Re}(\nu)$
belongs to the convex hull of
$\left\{w\rho\right\}_{\{w\in \rW\}}\subset \liea^{*}$,
a compact set.  Moreover, there exists $w\in\rW$ such that
$w^{2}=1$ and $w\nu=-\bar{\nu}$.  In particular if
$\mathrm{Re}(\nu)\neq 0$, then $w\neq 1$, and since
$\textrm{Im}(\nu)$ is $w$-fixed it is not regular.
\end{enumerate}
\end{thm}
Note that the norm on $\pi_{\nu}$ is only unique up to scaling. If
$\mathrm{Re}(\nu)=0$ and $\textrm{Im}(\nu)$ is regular (the main
case under consideration), we choose
$\left\Vert \varphi_{0}\right\Vert _{\pi_{\nu}}=1$.

For future reference we compute the action of $\mathfrak{g}$ on $V_{K}$
via $I_{\nu}$. First, remark that the action of $K$ on $V = \V$
is given by right translation, and the action of $\mathfrak{k} \subset
\mathfrak{g}$ on $V_K$ is then given by right differentiation.

Secondly, recall that if $U \subset \mathbb{R}^n$ is open,
a \emph{differential operator} $\mathbf{D}$ on $U$ is an
expression of the form
$\sum_{i=1}^{K} f_i \partial_1^{\alpha_1} \dots \partial_n^{\alpha_n}$,
where the $f_i$ are smooth and $\alpha_j \geq 0$. 
If $M$ is a smooth $n$-manifold,
we say a map $\mathbf{D}: C^{\infty}(M) \rightarrow C^{\infty}(M)$
is a differential operator if it is defined by a differential operator
in each coordinate chart. 

\begin{lem}
\label{lem: mathfrak{g}-on-ind}Let $f\in V_{K}$ and let $X\in\mathfrak{g}$.
Then there exists a differential operator $\mathbf{D}_X$ on $M \backslash K$
(depending linearly on $X$ and independent of $\nu$) such that for
every $k\in K$,\[
(I_{\nu}(X)f)(k)=\left\langle \nu+\rho,H_{0}(\Ad(k)X)\right\rangle f(k)+
(\mathbf{D}_X f)(k).\]
\end{lem}
\begin{proof}
Let $t\in\R$ be small, and consider $f(k\exp(tX))=f(\exp(t\Ad(k)X)\cdot k)$.
We write the Iwasawa decomposition of $\Ad(k)X\in\mathfrak{g}$ as
$\Ad(k)X=X_{\mathfrak{n}}(k)+X_{\mathfrak{a}}(k)+X_{\mathfrak{k}}(k)$
where $X_{\mathfrak{a}}(k)=H_{0}(\Ad(k)X)$. By the Baker-Campbell-Hausdorff
formula, $\exp(t\Ad(k)X)=\exp(tX_{\mathfrak{n}}(k))\cdot\exp(tX_{\mathfrak{a}}(k))\cdot\exp(tX_{\mathfrak{k}}(k))+O(t^{2})$,
so that:\[
(I_{\nu}(X)f)(k)=\frac{d}{dt}f\left(\exp(tX_{\mathfrak{n}}(k))\cdot\exp(tX_{\mathfrak{a}}(k))\cdot k\right)\restrict_{t=0}+\frac{d}{dt}f\left(\exp(tX_{\mathfrak{k}}(k))k\right)\restrict_{t=0}.\]
To conclude, observe that
$f \mapsto
\frac{d}{dt}f\left(\exp(tX_{\mathfrak{k}}(k))k\right)\restrict_{t=0}$
defines a differential operator $\mathbf{D}_X$ on $M \backslash K$. 
\end{proof}

Lemma \ref{lem: mathfrak{g}-on-ind} will be used in the following way:
as $\left\Vert \nu \right \Vert\rightarrow \infty$,
the operator 
$I_{\nu}(\frac{X}{\Vert \nu \Vert})$ acts on $V_K$
in a very simple fashion, \emph{modulo} certain error terms
of order $\Vert \nu \Vert^{-1}$. 
The simplicity of this ``rescaled''
action as $\Vert \nu \Vert \rightarrow \infty$ will be of importance in our
analysis. 

\subsection{\label{sub: functional} Some functional analysis} 
We collect here some simple functional analysis facts that we shall have
need of.

Let $\D$ denote the space of smooth functions of compact
support on $X$. It is endowed with the usual ``direct-limit'' topology:
fix a sequence of $K$-invariant compact sets $C_1 \subset C_2 \subset \dots$
such that their interiors exhaust $X$.  Then the $\DCi$ exhaust $\D$.
$\DCi$ is endowed as usual with a family of seminorms, viz. for any
$\mathcal{D}\in U(\lieg_\C)$ we define
$\left\Vert f\right\Vert_{C_i, \mathcal{D}} = \sup_{x\in C_i} |\mathcal{D}f|$.
These seminorms induce a topology on each $\DCi$. We give $\D$ the topology of
the union of $\DCi$, \ie a map from $\D$ is continuous if and only if its
restriction to each $\DCi$ is continuous.

In other words: a sequence of functions converges in $\D$
if their supports are all contained in a fixed compact set,
and all their derivatives converge uniformly on that compact set.

$\D$ is a locally convex complete space in this topology. In particular, its
subspace $\DK$ of $K$-finite vectors is dense.
We denote by $\D'$ (resp. $\DK'$) the topological dual to $\D$ (resp. the
algebraic dual to $\DK$).
Both spaces will be endowed with the weak-{*} topology.
We shall refer to an element of $\D'$ as a \emph{distribution} on $X$. 

Let $C_0(X)$ be the Banach space of continuous functions on $X$ decaying at
infinity, endowed with the supremum norm.  Let $C_0(X)'$ be the continuous
dual of $C_0(X)$; the Riesz representation theorem identifies it with the
space of finite (signed) Borel measures on $X$.
We endow $C_0(X)'$ with the weak-{*} topology. 

It is easy to see that $\DK$ is dense in $C_0(X)$. In particular any (algebraic)
linear functional on $\DK$ which is bounded \wrt the $\sup$-norm extends to
a finite signed measure on $X$, with total variation equal to the norm of the
functional. Moreover, if this functional is non-negative on the non-negative
members of $\DK$ then the associated measure is a positive measure.

\section{\label{sec: Lift}Representation-Theoretic Lift}

\subsection{\label{subsec:motivation}Introduction and motivation}
Suppose $\psi\in L^2(Y)$
has $\Vert\psi\Vert_2 = 1$ and an eigenfunction of $\mathfrak{Z}$. The
aim of the present section is to construct a distribution
$\mu_{\psi}$ on $Y$ that lifts the measure $\bar{\mu}_{\psi}$
on $Y$, and establish some basic properties of $\mu_{\psi}$. 

In the situation of Theorem \ref{mainthm}, if $\psi = \psi_n$,
the corresponding distribution will be the distributions
$\mu_n$ discussed in the proof of Theorem \ref{mainthm}.
The functions $\tilde{\psi}_n$ will then be chosen so that
the measures $|\tilde{\psi}_n(x)|^2 dx$ approximate $\mu_n$;
finally, both $|\tilde{\psi_n}(x)|^2$ and $\mu_n$ will become $A$-invariant
as $n \rightarrow \infty$.

We begin by fixing notation and providing
some motivation for the relatively formal definitions that follow. 

Setting $\psi(x)=\psi(xK)$ for any $x\in X$, we can think of $\psi$
as a function on $X$.  By the uniqueness of spherical functions
\cite[Th.\ 4.3 \& 4.5]{Helgason:GpsGeomAnal},
$\psi$ generates an irreducible spherical $G$-subrepresentation
of $L^{2}(X)$. As discussed in Section \ref{sub: repn-theory}, we
can then find $\nu\in\mathfrak{a}_{\C}^{*}$ such that this representation
is isomorphic to $\pi_{\nu}$ (in particular, $\pi_{\nu}$ is unitarizable).
We will assume for the rest of this section that $Re(\nu) = 0$,
\ie $\pi_{\nu}$ is \emph{tempered}, and that $\nu$ is regular. 
This will eventually be the only
case of interest to us in view of the non-degeneracy assumption
made later (Definition \ref{def: non-deg-simple}).
In this case $(V_{K},I_{\nu})$ is irreducible and isomorphic to $\pi_{\nu}$.
It follows that there is a unique $G$-homomorphism 
$R_{\psi}:(V,I_{\nu}) \to L^{2}(X)$ such that
$R_{\psi}(\varphi_{0})=\psi$.  The normalization $\left\Vert \psi\right\Vert _{L^{2}(X)}=1$
now implies $\left\Vert R_{\psi}(f)\right\Vert _{L^{2}(X)}=\left\Vert f\right\Vert _{L^{2}(K)}$
for any $f\in V_{K}$, \ie that $R_{\psi}$ is an isometry. 

We now give the rough idea of the construction that follows
in the language of Wolpert and Lindenstrauss;
the language we shall use later is slightly different,
so the discussion here also provides a translation.
The strategy of proof is similar to theirs; in a sense,
the main difficulty is finding the ``correct'' definitions in higher rank.
For instance, the proofs of Wolpert and Lindenstrauss use heavily
the fact that $K$-types for $\PSL_2(\R)$ have multiplicity one, 
and the explicit action of the Lie algebra by raising and lowering
operators.  We shall need a more intrinsic approach to handle the general
case.

The measure $\bar{\mu}_{\psi}$ on $Y$ is defined
by $g \mapsto \int_{X} g(x)|\psi(x)|^2 dx$.  More generally, suppose that
$\psi' \in L^2(X)$ belongs to the $G$-subrepresentation
generated by $\psi$, \ie $\psi' \in R_{\psi}(V)$.
We can then consider the (signed) measure 

\begin{equation} \label{eq:sigmadef}\sigma: g \mapsto \int_{X}
\psi(x) \overline{\psi'(x)} g(x) dx.\end{equation}
If $g(x)$ is $K$-invariant, then so is the product
$\psi(x) g(x)$, and it follows that the right-hand
side of (\ref{eq:sigmadef}) depends only on the projection
of $\psi'$ onto $R_{\psi}(V)^K$. 
The space $R_{\psi}(V)^K$ is one-dimensional, spanned by $\psi$,
and it follows that if $\psi' - \psi \perp \psi$ 
then the measure $\sigma$ on $X$ projects to the measure
$\bar{\mu}_{\psi}$ on $Y$.

The distribution $\mu_{\psi}$ we shall be construct
will be in the spirit of (\ref{eq:sigmadef}),
but with $\psi'$ a ``generalized vector'' 
in $R_{\psi}(V)$. Suppose, in fact, that $\psi'_1, \psi'_2, \dots
\psi'_n, \dots$ are an infinite sequence of elements of $R_{\psi}(V)$
that transform under different $K$-types, and suppose
further that $g \in \DK$. Then, by considering
$K$-types, the integral $\int_{X} \psi(x) \overline{\psi'_j(x)} g(x) dx$
vanishes for all sufficiently large $j$. It follows that, if one
sets $\psi'$ to be the \emph{formal sum} $\sum_{j=1}^{\infty} \psi'_j$,
one can make sense of (\ref{eq:sigmadef}) by interpreting it as:
$$\sigma(g) = \sum_{j=1}^{\infty} \int_{X} \psi(x) \overline{\psi'_j(x)} g(x)dx$$
In other words, if $g \in \DK$, we may make sense
of (\ref{eq:sigmadef}) while allowing $\psi'$ to belong to the space
$\widehat{V_K}$ of ``infinite formal sums of $K$-types.'' 
Our definition of $\mu_{\psi}$ will, indeed, be of the form
(\ref{eq:sigmadef}) but with $\psi'$ an ``infinite formal sum''
of this kind.  

For a certain choice of $\psi'$ (denoted $\Phi_\infty$ in
\cite{Lindenstrauss:HH_QUE}),
we will wish to
show that (\ref{eq:sigmadef}) is ``approximately a positive measure''
and ``approximately $A$-invariant,'' where both statements
become true in the large eigenvalue limit in an appropriate sense. 
For the ``approximate positivity,''
we shall integrate (\ref{eq:sigmadef}) by parts to show that
there exists another unit vector $\psi'' \in R_{\psi}(V)$ such that
$\int_{X} \psi(x) \overline{\psi'(x)} g(x) dx \approx 
\int_{X} |\psi''(x)|^2 g(x) dx$, where the right-hand
side is evidently a positive measure. For the ``approximate $A$-invariance,''
we will construct differential operators that annihilate
$\psi(x) \overline{\psi'(x)}$; this reduces to a purely algebraic
question of constructing elements in $U(\mathfrak{g})$ that
annihilate a vector in a certain tensor product representation. 

The space $\widehat{V_K}$ is very closely linked to the dual $V_K'$ of the
$K$-finite vectors: the conjugate linear isomorphism $T: V \rightarrow V'$
(\ref{eq:inclusion}) extends to a conjugate-linear isomorphism
$T: \widehat{V_K} \rightarrow V_K'$.  For formal reasons, it is simpler
to work with $V_K'$ than $\widehat{V_K}$; this is the viewpoint
we shall take in Definition \ref{def:mudef}. 
To motivate this viewpoint, let us rewrite (\ref{eq:sigmadef})
in a different fashion.
Let $v' \in V$ be chosen so that $\psi' = R_{\psi}(v')$, and let $P$
be the orthogonal projection of $L^2(X)$ onto $R_{\psi}(V)$.
We may rewrite (\ref{eq:sigmadef}) -- using the notations
of Definition \ref{def: vdef} -- as follows:
\begin{eqnarray} \nonumber \sigma(g)= \langle \psi(x) g(x), \psi'(x)
\rangle_{L^2(X)} =
\langle P(\psi(x) g(x)), \psi'(x) \rangle_{L^2(X)} \\ \label{eqn:froggy}
= \langle R_{\psi}^{-1} \circ P(\psi(x) g(x)), v' \rangle_{V} 
= T(v') \circ R_{\psi}^{-1} \circ P(\psi(x) g(x))\end{eqnarray}
Now, if $g \in \DK$, then the quantity
$R_{\psi} \circ P(\psi(x) g(x))$ is $K$-finite, \ie 
belongs to $V_K$. It follows
that, if $g \in \DK$, the last expression of (\ref{eqn:froggy}) makes formal
sense if we replace $T(v')$ by any functional $\Phi \in V_K'$.

\subsection{\label{sub: def-lift} Lifting a single (non-degenerate)
eigenfunction.}

\begin{defn} \label{def:mudef}
Let $\Phi \in V_K'$ be an (algebraic) functional, and $f \in V_{K}$.
Let $\mu_{\psi}(f,\Phi)$ be the functional on $\DK$
defined by the rule:
\begin{equation} \label{eq:mudef}
\mu_{\psi}(f, \Phi)(g) = \Phi 
\circ R_{\psi}^{-1} \circ P (R_{\psi} (f) \cdot g) \end{equation}
where $g \in \DK$, $P: L^2(X) \rightarrow R_{\psi}(V)$
is the orthogonal projection, and $R_{\psi}(f) \cdot g$
denotes pointwise multiplication of functions on $X$.
\end{defn}
\begin{rem}\label{rem: distr}
In fact, if $\Phi \in C^\infty(\mk)'$ (see equation (\ref{eq:inclusion})) then
$\mu_{\psi}(f,\Phi)$ extends to an element of $\D'$, \ie defines a
distribution on $X$: $\mu_\psi$ is the composite
$$\D \stackrel{g \mapsto R_{\psi}(f) g}{\longrightarrow} \D
\stackrel{R_{\psi}^{-1} P}{\longrightarrow} C^{\infty}(\mk)
\stackrel{\Phi}{\rightarrow} \C,$$
and it is easy to verify that each of these maps is continuous.
This is never used in our arguments: we use this observation only to
refer to certain $\mu_\psi$ as ``distributions''.
\end{rem}

\begin{defn}
\label{def: lift}Let $\delta\in V'_{K}$ be the distribution
$\delta(f)=f(1)$,
and call $\mu_{\psi} \stackrel{\mathrm{def}}{=}\mu_{\psi}(\varphi_{0},\delta)$ the (non-degenerate)
\emph{microlocal lift} of $\bar{\mu}_{\psi}$.
\end{defn}

The rest of the section will exhibit basic formal properties
of this definition. We will establish most of the formal
properties of $\mu_{\psi}$ by restricting $\Phi$ to be of the form
$T(f_2)$, where the conjugate-linear mapping $T$ is as defined
in (\ref{eq:inclusion}). This situation will occur sufficiently often that,
for typographical ease, it will be worth making the following definition:

\begin{defn}
Let $f_1, f_2 \in V_K$. We then set $\mu^{T}_{\psi}(f_1, f_2)
= \mu_{\psi}(f_1, T(f_2))$. 
\end{defn}

\begin{lem} \label{lem:mucty}Suppose $f_1, f_2 \in V_K$.  Then
\begin{equation}
\mu^{T}_{\psi}(f_{1},f_{2})(g)=\int_{X}R_{\psi}(f_{1})(x)
\overline{R_{\psi}(f_{2})(x)} g(x)dx.\label{eq: def-mu_R}\end{equation}
and $\mu^{T}_{\psi}$ defines a signed measure on $X$ of total variation
at most $||f_1||_{L^2(K)} ||f_2||_{L^2(K)}$. 
If $f_1 = f_2$, then $\mu^{T}_{\psi}(f_1, f_1)$ is a positive 
measure of mass $||f_1||^2_{L^2(K)}$. 
\end{lem}
\proof  (\ref{eq: def-mu_R}) is a consequence of the definition of $\mu$.
The Cauchy-Schwarz inequality implies that $|\mu^T_{\psi}(f_1, f_2)(g)|
\leq ||f_1||_{L^2(K)} ||f_2||_{L^2(K)} ||g||_{L^{\infty}(X)}$, whence
the second conclusion. The last assertion is immediate. 
\qed

In fact, it may be helpful to think of $\mu_{\psi}$
as being given by a distributional extension
of the formula (\ref{eq: def-mu_R}); see the discussion 
of Section \ref{subsec:motivation}. 

\begin{lem} \label{lem:projects}
The distribution $\mu_{\psi}(\varphi_0, \delta)$ on $X$
projects to the measure $|\psi|^2 dy$ on $Y$.
\end{lem}
\begin{proof} 
In view of the previous Lemma, it will suffice
to show that the distribution $\mu_{\psi}(\varphi_0, \delta)
- \mu_{\psi}^T(\varphi_0, \varphi_0)$ on $X$ projects to 
$0$ on $Y$. This amounts to showing that $\mu_{\psi}(\varphi_0,
\delta - T(\varphi_0))$ annihilates any $K$-invariant
function $g \in \D^K$.  Taking into account
that the functional $\delta - T(\varphi_0)$ on $V_K$
annihilates any $K$-invariant vector, 
the claim follows
from the definition of $\mu_{\psi}$. 
\end{proof}

\begin{lem}\label{lem:equiv}
The map $\mu_{\psi}: V_K \otimes V_K' \rightarrow \DK'$
is equivariant for the natural $\mathfrak{g}$-actions on both sides. 
\end{lem}
\begin{proof} This follows directly from the definition of $\mu_{\psi}$. 
\end{proof}

Concretely speaking, this says that for $f \in V_K, \Phi \in V_K',
g \in \DK, X \in \mathfrak{g}$ we have
\begin{equation} \label{eq:concrete}\mu_{\psi}(Xf_{1},\Phi)(g)+
\mu_{\psi}(f_{1},X\Phi)(g)+\mu_{\psi}(f_{1},\Phi)(Xg)=0\end{equation}
where $X$ acts on $V_{K}$ via $I_{\nu}$ and on $V_K'$ via
$I_{\nu}'$. In particular, if $f_1, f_2 \in V_K$ we have
\begin{equation} \label{eq:concrete2}\mu^T_{\psi}(Xf_{1},f_2)(g)+
\mu_{\psi}^T(f_{1},Xf_2)(g)+\mu_{\psi}^T(f_{1},f_2)(Xg)=0\end{equation}

\subsection{\label{sub: positivity}Sequences of eigenfunctions and quantum limits.}
In what follows we shall consider
$\left\{ \psi_{n}\right\} _{n=1}^{\infty}\subset L^{2}(Y)$, 
a sequence of eigenfunctions
with parameters $\left\{ \nu_{n}\right\} $ diverging to $\infty$
(\ie leaving any compact set). Set $\tilde{\nu}_n
= \frac{\nu_n}{||\nu_n||}$. For $f_1, f_2 \in V_K$ and $\Phi \in V_K'$,
we abbreviate $\mu^T_{\psi_n}(f_1, f_2)$ 
(resp. $\mu_{\psi_n}(f ,\Phi)$) to $\mu^T_n(f_1, f_2)$ 
(resp. $\mu_n(f, \Phi)$), and 
we abbreviate the microlocal lift $\mu_{\psi_n}$ ($:= \mu_n(\varphi_0,
\delta)$) to $\mu_n$.

\begin{defn}
\label{def: non-deg-simple}($G$ simple) 
We say a sequence $\psi_n$ is \emph{non-degenerate}
if every limit point of the sequence $\tilde{\nu}_n$ is regular.

We say that it is \emph{conveniently arranged}
if it is nondegenerate, $\lim_{n \rightarrow \infty} \tilde{\nu}_n$
exists, $Re(\nu_n) = 0$ for all $n$,
the $\nu_n$ are all regular, and for all
$f_1, f_2 \in V_K$ the measures $\mu_n^T(f_1, f_2)$ converge
in $C_0(X)'$ as $n \rightarrow \infty$.  
In this situation we denote $\lim_{n \rightarrow \infty} \tilde{\nu}_n$
by $\tilde{\nu}_{\infty}$. 
\end{defn}

The existence of non-degenerate sequences of eigenfunctions was discussed
in Remark \ref{rem: mainthm}.  This follows from strong versions of Weyl's
Law on $Y$.  By Theorem \ref{thm: ss-gps-theory},
the non-degeneracy of a sequence $\psi_n$ as in the Definition
implies $\mathrm{Re}(\nu_{n})=0$
for all large enough $n$. 
For fixed $f_1, f_2 \in V_K$ the total variation of the measures 
$\mu_n^{T}(f_1, f_2)$ is bounded independently of $n$ (Lemma
\ref{lem:mucty});
in view of the (weak-{*}) compactness of the 
unit ball in $C_0(X)'$ it follows that this sequence
of measures has a convergent subsequence.
Combining this remark with the fact that $V_K$ has a countable basis,
a diagonal argument shows that every non-degenerate sequence of eigenfunctions has a conveniently arranged
subsequence. 

Now suppose $\{\psi_n\}$ is a conveniently arranged sequence and fix
$f_1 \in V_K, \Phi \in V_K', g \in \DK$. Let $\Phi_N$ be 
the $N$-truncation of $\Phi$ (see Definition \ref{def: vdef}). 
In view of (\ref{eq:mudef}), if we  choose $N := N(f_1, g)$ sufficiently
large, then $\mu_{n}(f_1, \Phi)(g) = \mu^T_{n}(f_1, \Phi_N)(g)$.
It follows that the limit $\lim_{n \rightarrow \infty} \mu_n(f_1, \Phi)(g)$
exists.

We may consequently define
$\mu_{\infty}\colon V_K \times V_K' \to \DK'$
and $\mu_{\infty}^T\colon V_K \times V_K \to \DK'$ by the rules:
\begin{equation} \label{eqn:muinftydef}
\begin{aligned}
\mu_{\infty}(f, \Phi)(g) &= \lim_{n\to\infty} \mu_n(f_1, \Phi)(g),
\quad (g \in \DK) \\
\mu^T_{\infty}(f_1, f_2) &= \mu_{\infty}(f_1, T(f_2)) &
\end{aligned}
\end{equation}

\begin{lem} \label{lem:mucty2}
For fixed $f_1 \in V_K$, the map $\Phi \rightarrow 
\mu_{\infty}(f_1, \Phi)$ is continuous as a map
$V_K' \rightarrow \DK'$, both spaces being endowed with the weak
topology. 
\end{lem}
\begin{proof} This is an easy consequence of the definitions.
\end{proof}
It is natural to ask whether $\mu_{\infty}(f_1, \Phi)$
extends to an element of $\D'$, especially if $\Phi \in C^{\infty}(\mk)'$.
Indeed, it is possible to make quantitative the argument of
Remark \ref{rem: distr} to obtain a uniform bound on the
distributions $\mu_n(f_1,\Phi)$.  This will not be needed in
this paper, however, since for our choice of $(f_1,\Phi)$,
the limiting distribution is positive (in particular a measure), a fact
we will prove directly.

Henceforth $\left\{ \psi_{n}\right\} _{n=1}^{\infty}$
will be a conveniently arranged sequence.
We will show that $\mu_{\infty}(\varphi_{0},\delta)$
is positive and bounded \wrt the $L^{\infty}$ norm on $\DK$.
It hence extends to a finite positive measure.

\begin{rem*}
\textbf{~}
In the case of a semisimple group, one can allow the projection of
the parameter to each simple factor to tend to infinity at a different
rate. The definition of a non-degenerate limit can then remain unchanged.
The $\tilde{\nu}_{n}$ however must be defined with greater care --
see Section \ref{sub: extend}. 
\end{rem*}
The key to the positivity of the limits is the following lemma (cf. 
\cite[Prop.\ 3.3]{Wolpert:SL2_Lift_Fejer},
\cite[Th.\ 3.1]{Lindenstrauss:HH_QUE}).

\begin{lem}
\label{lem: int-parts}(Integration by parts) Let $\left\{ \psi_{n}\right\} $
be conveniently arranged. Then, for any $f,\, f_{1},\, f_{2}\in V_{K}$ we
have:
\begin{equation} \label{eq:byparts} 
\mu^T_{\infty}(f_{1},  f \cdot f_{2})
=\mu^T_{\infty}(\bar{f} \cdot f_{1},f_{2}).\end{equation}
Here \eg $f \cdot f_2$ denotes pointwise multiplication of functions on $M
\backslash K$. 
\end{lem}
\begin{proof}
We start by exhibiting explicit
functions $f$ for which (\ref{eq:byparts}) is valid.

Extend every $\nu\in\mathfrak{a}_{\C}^{*}$ to $\mathfrak{g}_{\C}^{*}$
via the Iwasawa decomposition $\mathfrak{g}=\mathfrak{n}\oplus\mathfrak{a}\oplus\mathfrak{k}$.
For any $X\in\mathfrak{g}$, let $p_{X}(k)=\frac{1}{i}\left\langle \tilde{\nu}_{\infty},\Ad(k)X\right\rangle $.
For fixed $X$, $k \mapsto p_{X}(k)$ defines a $K$-finite element of $\V$.

By (\ref{eq:concrete2}), for every $X$, $f_{1},$ $f_{2}$,
$g$, and $n$, we have 
\begin{equation} \mu_{n}^T(Xf_{1},f_{2})(g)+
\mu_{n}^T(f_{1},X \,  f_{2})(g)+\mu_{n}^T(f_{1},f_{2})(Xg)=0. \end{equation}
Divide by $\left\Vert \nu_{n}\right\Vert $ and apply
Lemma \ref{lem: mathfrak{g}-on-ind} to see:
\begin{eqnarray} \label{eq:verylong}
\mu^T_{n}(ip_{n}\cdot f_{1},f_{2})(g)+\mu^T_{n}(f_{1},ip_{n} \cdot f_{2})(g)
\\ \nonumber
=-\frac{ \mu^T_{n}(\mathbf{D}_X f_{1},f_{2})(g)+
\mu_{n}^T(f_{1},\mathbf{D}_X f_{2})(g)+
\mu_{n}^T(f_{1},f_{2})(X g)}{\left\Vert
\nu_{n}\right\Vert },\end{eqnarray}
where $p_{n}(k)=\frac{1}{i}\left\langle \tilde{\nu}_{n}+\frac{\rho}{\left\Vert \nu_{n}\right\Vert },\Ad(k)X\right\rangle $.

As $n\to\infty$, the right-hand side of (\ref{eq:verylong}) tends to
zero by Lemma \ref{lem:mucty}. 
On the other hand,  $p_{n}f_{i}$ (considered as continuous functions on $K$)
converge uniformly to $p_{X}f_{i}$.
Another application of Lemma \ref{lem:mucty} shows that
the left-hand side of (\ref{eq:verylong}) converges to
$i\mu_{\infty}^T(p_{X}f_{1}, f_{2})-i\mu_{\infty}^T(f_{1},p_{X} \cdot f_{2})$.
Since $p_X = \overline{p_X}$ this
shows that (\ref{eq:byparts}) holds with $f = p_{X}$.  

Now let $\mathcal{F}\subset C(\mk)$ be the $\C$-subalgebra generated
by the $p_{X}$ and the constant function $1$. Clearly
(\ref{eq:byparts}) holds for all $f\in\mathcal{F}$. This subalgebra is
$K$-stable since $p_{X}(k k_1) = p_{\mathrm{Ad}(k_1) X}(k)$ and hence $\mathcal{F}\cap V_{\rho}\subset\mathcal{F}$
for all $\rho\in\hat{K}$. Showing $\mathcal{F}$ is dense in $\V$
suffices to conclude that $\mathcal{F}=V_{K}$. 

We will prove the
stronger assertion that $\mathcal{F}$ is dense in $C(\mk)$ using
the Stone-Weierstrass theorem. 
Note that $1\in\mathcal{F}$, and $\mathcal{F}$ is closed under complex
conjugation since $p_X = \overline{p_X}$. 
It therefore suffices to show that
$\mathcal{F}$ separates the points of $\mk$. To this end, let $k_{1},k_{2}\in K$
be such that $p_{X}(k_{1})=p_{X}(k_{2})$ for all $X\in\mathfrak{g}$.
Then $\left\langle \tilde{\nu}_{\infty},\Ad(k_{1})X\right\rangle =\left\langle \tilde{\nu}_{\infty},\Ad(k_{2})X\right\rangle $
for all $X\in\mathfrak{g}$,
\ie $\langle \Ad(k_1)^{-1} \tilde{\nu}_{\infty} - Ad(k_2)^{-1}
\tilde{\nu}_{\infty}, X \rangle = 0$ for all $X \in \mathfrak{g}$.
This implies that $\Ad(k_1^{-1}) \tilde{\nu}_{\infty} = \Ad(k_2)^{-1}
\tilde{\nu}_{\infty}$; by the non-degeneracy assumption,
$Z_K(\tilde{\nu}_{\infty}) = Z_K(A) = M$, so $Mk_1 = M k_2$,
\ie $k_1$ and $k_2$ represent the same point of $\mk$.
\end{proof}

Lemma \ref{lem: int-parts} shows easily that
$\mu_{\infty}(\varphi_0, \delta)$ extends to a positive measure. 
Indeed, choosing $f_j$ as in Lemma \ref{lem: sqrt-delta},
we see that 
\begin{equation} \label{eq:baby}\mu_{\infty}(\varphi_0, \delta) = \lim_{j \to \infty} \mu^T_{\infty}(\varphi_0,
|f_j|^2) = \lim_{j \to \infty} \mu_{\infty}^T(f_j, f_j).\end{equation}
Here we have invoked Lemma \ref{lem:mucty2} for the first equality.
It is clear that $\mu_{\infty}^T(f_j, f_j)$ defines a positive measure on
$X$; thus $\mu_{\infty}(\varphi_0, \delta)$, initially
defined as an (algebraic) functional on $\DK$, extends to a positive measure on
$X$. 
To obtain the slightly stronger conclusion implicit
in (\ref{thm:claim2}) of Theorem \ref{mainthm}, we will analyze 
this argument more closely.

\begin{cor}
\label{cor: convergence}
Notations as in Lemma \ref{lem: int-parts},
there exist a constant $C_{f_1,f_2,f}$
and a seminorm $||\cdot||$ on $\D$ such that
\begin{equation} \label{eq:hard}
\left| \mu_n^T(f_1, f \cdot f_2)(g) - \mu_n^T(\overline{f} \cdot f_1,
f_2)(g) \right| \leq C_{f_1, f_2, f} \Vert g\Vert \left[\Vert\tilde{\nu}_{\infty}
- \tilde{\nu}_n\Vert + \Vert\nu_n\Vert^{-1}\right]\end{equation}
\end{cor}
\begin{proof}
This follows by 
keeping track of the error term in the proof of
of Lemma \ref{lem: int-parts}. 

Fix a basis $\{X_i\}$ for $\lieg$, and define
a seminorm on $\D$ by $\Vert g \Vert =
\Vert g \Vert_{L^\infty(X)} + \sum_{i} \Vert X_i g \Vert_{L^\infty(X)}$. 
With this seminorm, (\ref{eq:hard}) holds for $f_1, f_2 \in V_K$
and $f = p_X$.  This follows from (\ref{eq:verylong}),
utilizing Lemma \ref{lem:mucty} and the fact that
$\Vert p_X - p_n \Vert_{L^\infty(\mk)} \ll
\Vert\tilde{\nu}_\infty - \tilde{\nu}_n \Vert$.

Next suppose $f_1, f_2, f, f' \in V_K$ and
$\alpha, \alpha' \in \C$. Then, if (\ref{eq:hard})
is valid for $(f_1, f_2, f)$ and $(f_1, f_2, f')$,
it is also valid for $(f_1, f_2, \alpha f + \alpha' f')$.
Further, if (\ref{eq:hard}) is valid for $(f_1, f' \cdot f_2, f)$
and for $(\overline{f}  f_1, f_2, f')$, then it is also valid
for $(f_1, f_2, f \cdot f')$.

Consider now the set of $f \in V_K$ for which (\ref{eq:hard})
holds for all $f_1, f_2 \in V_K$. The remarks above
show that this is a subalgebra of $V_K$ that contains each
$p_X$. The Corollary then
follows from the equality $\mathcal{F}=\V_K$ established in the Lemma.
\end{proof}

\begin{rem} \label{rem:subseq}
In fact, it is possible to obtain a bound of the form
$C_{f_1,f_2,f,\tilde{\nu}_n}\Vert g\Vert \Vert\nu_n\Vert^{-1}$, with the
constant uniformly bounded if the $\tilde{\nu}_n$ are uniformly bounded away
from the walls.  This result can be used
to avoid passing to a subsequence in Theorem \ref{mainthm}
or the following Proposition; this is unnecessary for our applications,
however. 
\end{rem}

\begin{prop}
\label{prop: equivariance}(Positivity
and equivariance:  (\ref{thm:claim2}) and (\ref{prop:eq})
of Theorem \ref{mainthm}).

Let $\left\{ \psi_{n}\right\} $
be non-degenerate. 
After replacing $\left\{ \psi_n \right\}$ by an appropriate subsequence,
there exist functions $\tilde{\psi}_n$ on $X$ with the following properties:
\begin{enumerate}
\item  \label{claim1} Define the measure $\sigma_n$
via the rule $\sigma_n(g) = 
\int_{X} g(x) |\tilde{\psi_n}(x)|^2 dx$.
Then, for each $g \in \DK$ we have
$\lim_{n \rightarrow \infty} (\sigma_n(g) - \mu_n(g))  = 0$.

\item \label{claim3}
Let $E \subset \mathrm{End}_G(\CX)$
be a $\C$-subalgebra of endomorphisms of $\CX$, commuting
with the $G$-action.  Note that each $e \in E$ induces an endomorphism
of $C^{\infty}(Y)$.
Assume in addition that $\psi_n$ is an eigenfunction for $E$.
Then we may choose $\tilde{\psi}_n$ so that each
$\tilde{\psi}_n$ is an eigenfunction for $E$
with the same eigenvalues as $\psi_n$.
\end{enumerate}
\end{prop}

\begin{proof}
Without loss of generality we may assume that $\{\psi_n \}$ are
conveniently arranged. 

Let $\left \{ f_{j}\right\} _{j=1}^{\infty}\subset V_{K}$
be the sequence of functions provided by Lemma \ref{lem: sqrt-delta},
so that $T(|f_j|^2)$ approximates $\delta$.  The main idea is,
as in (\ref{eq:baby}), to
approximate $\mu_n = \mu_n(\varphi_0,\delta)$
using $\mu^T_n(f_j,f_j)$.

For any $g\in\DK$ we have:
\begin{eqnarray} \label{eq:mubound}
\left|\mu_{n}(g) - \mu^T_n(f_j, f_j)(g)\right| \leq
 \left|\mu_{n}(\varphi_0, \delta)(g) - \mu_n(\varphi_0, |f_j|^2)(g)\right|
\nonumber \\ + \left|\mu_n(\varphi_0, |f_j|^2)(g) - \mu_n(f_j, f_j)(g)\right|.
\end{eqnarray}

Corollary (\ref{cor: convergence}) provides a seminorm
$\left\Vert\cdot\right\Vert$ on $\D$ and a constant $C_j$ such that
$$\left|\mu_n(\varphi_0, |f_j|^2)(g) - \mu_n(f_j, f_j)(g)\right| \leq 
C_j ||g|| \cdot \left[
\left\Vert\tilde{\nu}_n - \tilde{\nu}_{\infty}\right\Vert +
\left\Vert\nu_n\right\Vert^{-1} \right].$$
Choose a sequence of integers $\{j_n\}_{n=1}^\infty$ such that
$j_n \rightarrow \infty$ and:
$$C_{j_n} \cdot \left[
\left\Vert\tilde{\nu}_n - \tilde{\nu}_{\infty}\right\Vert +
\left\Vert\nu_n\right\Vert^{-1} \right]
\xrightarrow[n\to\infty]{} 0$$

We now estimate the other term on the right-hand side of (\ref{eq:mubound}).
Choosing $N = N(g)$ large enough so that
$\mu_n(\varphi_0,\delta)(g) = \mu_n(\varphi_0,\delta_N)(g)$, we have
$$\left|\mu_{n}(\varphi_0, \delta)(g) - \mu_n(\varphi_0, |f_j|^2)(g)\right| \leq \left\Vert |f_j|^2_N - \delta_N \right\Vert_{\V}
\left\Vert g \right\Vert_\infty.$$  As $j\to\infty$ (in particular, if
$j=j_n$), $|f_j|^2_N \to \delta_N$ in $V_N$, so this term tends to zero.
It follows that \begin{equation}
\label{eq:strom1} \lim_{n\to\infty}
\left|\mu_n(g) - \mu^T_n(f_{j_n},f_{j_n})(g)\right| = 0,
\end{equation}

Setting $\tilde{\psi}_n = R_{\psi_n}(f_{j_n})$, we deduce that
\begin{equation} \label{eq:strom3} 
\lim_{n \to \infty} \left(\mu_{n}(g) -
\int_{X} |\tilde{\psi}_n|^2 g(x) dx \right) = 0 \end{equation}
holds for every $g \in \DK$.  In particular, we obtain (\ref{claim1})
of the Proposition.

To obtain the equivariance property note that the representation
$I_{\nu_n}$ is irreducible as a $\left(\lieg,K\right)$-module.
By \cite[Corollary 8.11]{Knapp:Rpn_Th_SS_Gps}, there
exists $u_n \in U(\lieg)$ such that $I_{\nu_n}(u_n) \varphi_0 = f_{j_n}$.
Thus $\tilde{\psi}_n = u_n \psi_n$. Now every $e \in E$ commutes
with the right $G$-action; in particular, $e u_n = u_n e$.
It follows that $\tilde{\psi}_n$ transforms under the same character
of $E$ as $\psi_n$.
\end{proof}

\section{\label{sub: A-inv}Cartan invariance of quantum limits}

In this section we show that a nondegenerate quantum limit $\mu_{\infty}$
is invariant under the action of $A<G$.
This invariance follows from differential equations
satisfied by the intermediate distributions $\mu_{n}$. The construction
of these differential equations is a purely algebraic problem:
construct elements in the $U(\lieg_{\C})$-annihilator of $\varphi_0 
\otimes \delta \in V_K \otimes V_K'$, where the $U(\lieg_{\C})$-action
is by $I_{\nu} \otimes I_{\nu}'$.

Ultimately, these differential equations are derived from the fact that each
$z\in\mathfrak{Z}=\mathfrak{Z}(\mathfrak{g}_{\C})$
acts by a scalar on the representation $(V_{K},I_{\nu_{n}})$. 
To motivate the method and provide an example, we
first work out the simplest case, that of $\PSL_{2}(\R)$, in detail.  In this
case the resulting operator is due to Zelditch.

\subsection{Example of $G=\PSL_{2}(\R)$.}

Set $G = \PSL_2(\R)$, $\Gamma \leq G$ a lattice,
and $A$ the subgroup of diagonal matrices. 
Let $H$ (explicitly given below) be the infinitesimal generator of
$A$, thought of first as a differential operator acting on $X=\Gamma\backslash G$
via the differential of the regular representation. 
If $\{\psi_n\}$ is a conveniently arranged sequence of eigenfunctions on 
$\Gamma \backslash G/K$, and $\mu_n$ the corresponding
distributions (Definition \ref{def: lift}),
we will exhibit
a second-order differential operator $J$ such that for all $g\in\DK$,
\begin{equation}
\mu_{n}((H-\frac{J}{r_{n}})g)=0,\label{eq: zelditch-op}\end{equation}
 where $r_{n}\sim|\lambda_{n}|^{1/2}$. Since the $\mu_{n}(Jg)$ are
bounded (they converge to $\mu_{\infty}(Jg)$), we will conclude that
$\mu_{\infty}(Hg)=0$, in other words that $\mu_{\infty}$ is $A$-invariant.
This operator in equation (\ref{eq: zelditch-op}) is given in
\cite{Zelditch:SL2_Lift_A_inv}. Its discovery was motivated by the proof
(via Egorov's theorem) of the invariance
of the usual microlocal lift under the geodesic flow. We show here
how it arises naturally in the representation-theoretic approach.

By Lemma \ref{lem:equiv}, it will suffice to find an operator
annihilating the element $\varphi_{0}\otimes\delta \in V_{K}\otimes 
V'_{K}$,
where $U(\mathfrak{g}_{\C})$ acts via $I_{\nu}\otimes I_{\nu}'$.

Let $H=\left(\begin{array}{cc}
1\\
 & -1\end{array}\right)$, $X_{+}=\left(\begin{array}{cc}
0 & 1\\
 & 0\end{array}\right)$, $X_{-}=\left(\begin{array}{cc}
0\\
1 & 0\end{array}\right)$ be the standard generators of $\mathfrak{sl}_{2}$, with the commutation relations
$\left[H,X_{\pm}\right]=\pm2X_{\pm}$, $\left[X_{+},X_{-}\right]=H$.
The roots \wrt the maximal split torus $\mathfrak{a}=\R\cdot H$
are given by $\pm\alpha(H)=\pm2$.
We also set $W=X_{+}-X_{-}$, so that $\R\cdot W=\mathfrak{k}$.
Letting $+\alpha$ be the positive root, $\mathfrak{n}=\R\cdot X_{+}$,
we have $\rho(H)=\frac{1}{2}\alpha(H)=1$.
Set $\exp\mathfrak{a}=A$ as in
the introduction. 

The Casimir element $C\in\mathfrak{Z}(\mathfrak{sl}_{2}\C)$ is given
by $4C=H^{2}+2X_{+}X_{-}+2X_{-}X_{+}$. For the parameter $\nu\in i\mathfrak{a}^{*}$
given by $\nu(H)=2ir$ ($r\in\R$), $C$ acts on $\pi_{\nu}$ with
the eigenvalue $\lambda=-\frac{1}{4}-r^{2}$. The Weyl element acts
by mapping $\nu\mapsto-\nu$.  On $S=G/K$
with the metric normalized to have constant curvature $-1$, $C$
reduces to the hyperbolic Laplacian. In particular, every eigenfunction
$\psi\in L^{2}(\Gamma\backslash G/K)$ with eigenvalue $\lambda<-\frac{1}{4}$
generates a unitary principal series subrepresentation. 
Definition \ref{def: lift} associates
to $\psi$ a distribution $\mu_{\psi}(\varphi_{0},\delta)$ on $\Gamma
\backslash G$.

As in Definition \ref{defn:repdefn}, we have an action
$I_{\nu}$ of $G$ on $V$ and of $\mathfrak{g}$ on $V_K$. 
Note that for $g\in NA$, $f\in V_K$,
$\left(I_{\nu}(g)f\right)(1)=f(g)=e^{\left\langle
\nu+\rho,H_{0}(g)\right\rangle }f(1)$.
Since $ \delta(f) =f(1)$ and the pairing between $V_K$ and $V_K'$
is $G$-invariant, it follows that for $X\in\mathfrak{a}\oplus\mathfrak{n}$,
$I_{\nu}'(X)\delta=-\left\langle \nu+\rho,H_{0}(X)\right\rangle \delta$.

Suppressing $I_{\nu}$ from now on this means that 
$X\cdot(f\otimes\delta)=(Xf)\otimes\delta-\left\langle
\nu+\rho,H_{0}(X)\right\rangle f\otimes\delta$.
Extend $\nu+\rho$ trivially on $\mathfrak{n}$ to obtain
a functional on $\mathfrak{a}\oplus\mathfrak{n}$. 
Then 
\begin{equation} \label{eq:one}
(X+(\nu+\rho)(X))\cdot (f\otimes\delta) = (Xf)\otimes\delta.
\end{equation}
Now since $\mathfrak{a}$ normalizes $\mathfrak{n}$ and $\nu+\rho$
is trivial on $\mathfrak{n}$, the map $X\mapsto X+\left(\nu+\rho\right)(X)$
is a Lie algebra homomorphism $\mathfrak{a}\oplus\mathfrak{n}\to\mathfrak{a}\oplus\mathfrak{n}$,
and hence extends to an algebra homomorphism $\tau_{\nu+\rho}\colon U(\mathfrak{a}_{\C}\oplus\mathfrak{n}_{\C})\to U(\mathfrak{a}_{\C}\oplus\mathfrak{n}_{\C})$.
(\ref{eq:one}) shows that, for $u \in U(\mathfrak{a}_{\C} \oplus
\mathfrak{n}_{\C})$,
\begin{equation}
\label{eq:tens}\tau_{\nu+\rho}(u)\cdot( f
\otimes\delta)=(uf)\otimes\delta\end{equation}
In view of (\ref{eq:tens}) any operator $u \in U(\mathfrak{a}_{\C}
\oplus \mathfrak{n}_{\C})$ annihilating $\varphi_0$ gives
rise to an operator annihilating $\varphi_{0} \otimes \delta$. 

The natural starting point is the eigenvalue equation
$(4C+1+4r^{2})\varphi_{0}=0$. Of course, $C$ is not
an element of $U(\mathfrak{n}_{\C}\oplus\mathfrak{a}_{\C})$. Fortunately,
it ``nearly'' is: there exists an $C' \in U(\mathfrak{n}_{\C} 
\oplus \mathfrak{a}_{\C})$ such that $C - C'$ annihilates $\varphi_0$. 

In detail, we use the commutation relations and the fact that $X_{-}=X_{+}-W$
to write $4 C = H^{2}-2H+4X_{+}^{2}-4X_{+}W$. 
Since $\phi_{0}$ is spherical, it follows that  $W \phi_0 = 0$.  
Thus 
\begin{equation} \label{eq:almostthere}
\left(H^{2}-2H+4X_{+}^{2}+1+4r^{2}\right)\varphi_{0}=0\end{equation}
Since $(\nu+\rho)(H)=2ir+1$, we conclude from (\ref{eq:tens})
that:\[
\left((H+2ir+1)^{2}-2(H+2ir+1)+4X_{+}^{2}+1+4r^{2}\right)\cdot
\varphi_{0}\otimes \delta=0.\]
Collecting terms in powers of $r$ we see that this may be written as:
\[
\left((2H)(2ir)+(H^{2}+4X_{+}^{2})\right) \varphi_{0} \otimes \delta  = 0\]
Setting $J=\frac{H^{2}+4X_{+}^{2}}{4i}$ and dividing by $4ir$ we
see that the operator $H + \frac{J}{r}$ annihilates $\varphi_{0}
\otimes \delta$,  and so also the distribution $\mu_n$. 
One then deduces the $A$-invariance of $\mu_{\infty}$ as discussed
in the start of this section. 

Notice that the terms involving
$r^2$ in (\ref{eq:almostthere}) canceled. This is a general feature
which will be of importance.

\subsection{The general proof}

We now generalize these steps in order. Notations
being as in Section \ref{sec: Notation} and 
Definition \ref{def: lift}, 
we first compute the action
of $U(\mathfrak{m_{\C}}\oplus\mathfrak{a}_{\C}\oplus\mathfrak{n}_{\C})$
on $\delta$ (Lemma \ref{lem: U-on-phi_inf}) and then on $\varphi_{0}\otimes
\delta$
(Corollary \ref{cor: U-on-phi-phi_inf}). Secondly we find an appropriate
form for the elements of $\mathfrak{Z}(\mathfrak{g}_{\C})$ (Corollary
\ref{cor: z-pr(z) rewrite}), which gives us the exact differential
equation (\ref{eq: ann-phi_0-phi_inf}). We then show that the elements
we constructed annihilating $\mu_{\psi}$ are (up to scaling) of an
appropriate form $H+\frac{J}{\left\Vert \nu\right\Vert ^{*}}$ (Lemma
\ref{lem: calc-J}), and ``take the limit as $\nu\rightarrow\infty$''
(Corollary \ref{cor: P'-inv}) to see that $\mu_{\infty}$ is invariant
under a sub-torus of $A$. 

A final step (not so apparent in the $\PSL_2(\R)$ case)
is to verify that 
we have constructed \emph{enough} differential operators to obtain
invariance under the full split torus (Lemma \ref{lem:previous}).
In fact, even in the rank-$1$ case one needs to verify that the ``$H$''
part is non-zero.

Given $\lambda\in\mathfrak{a}_{\C}^{*}$, we extend it to a linear
map $\mathfrak{m}_{\C}\oplus\mathfrak{a}_{\C}\oplus\mathfrak{n}_{\C}\to\C$.
Since $\mathfrak{m}_{\C}\oplus\mathfrak{n}_{\C}$ is an ideal of this
Lie algebra, $\lambda$ is a Lie algebra homomorphism;
thus it extends to an algebra homomorphism $\lambda\colon U(\mathfrak{m}_{\C}\oplus\mathfrak{a}_{\C}\oplus\mathfrak{n}_{\C})\rightarrow\C$.
We denote by $\tau_{\lambda}$ the translation automorphism of $U(\mathfrak{m}_{\C}\oplus\mathfrak{a}_{\C}\oplus\mathfrak{n}_{\C})$
given by $X\mapsto X+\lambda(X)$ on $\mathfrak{m}_{\C}\oplus\mathfrak{a}_{\C}\oplus\mathfrak{n}_{\C}$.
Similarly, given $\chi\in\mathfrak{h}_{\C}^{*}$, we define $\tau_{\chi}:U(\mathfrak{h}_{\C})\rightarrow U(\mathfrak{h}_{\C})$.
We shall write $U(\mathfrak{g}_{\C})^{\leq d}$ for the elements of
$U(\mathfrak{g}_{\C})$ of degree $\leq d$, and similarly for other
enveloping algebras and $\mathfrak{Z}=\mathfrak{Z}(\mathfrak{g}_{\C})$
(\eg $\mathfrak{Z}^{\leq d}=\mathfrak{Z}\cap U(\mathfrak{g}_{\C})^{\leq d}$).

Let $\nu \in \mathfrak{a}_{\C}^{*}$. 
Let $\chi_{\nu}:\mathfrak{Z}\rightarrow\C$ be the infinitesimal character
corresponding to $I_{\nu}$ (that is, the scalar by which $\mathfrak{Z}$
acts in $(I_{\nu},V_K)$.)
Recall that $\rho_{\mathfrak{h}}$ denotes
the half-sum of positive roots for $\hg$, $\rho$ the half-sum for
$\ag$.

\begin{lem}
\label{lem: U-on-phi_inf}For $X\in\mathfrak{m}\oplus\mathfrak{a}\oplus\mathfrak{n}$,
$I_{\nu}(X)\delta=-\left\langle \nu+\rho,X\right\rangle \delta$.
\end{lem}
\begin{proof} This follows from the definitions.
\end{proof}

\begin{cor}
\label{cor: U-on-phi-phi_inf}For any $u \in U(\mathfrak{m}_{\C}\oplus\mathfrak{a}_{\C}\oplus\mathfrak{n}_{\C})$
and $f\in V_{K}$,
$$I_{\nu} \otimes
I_{\nu}'(\tau_{\nu+\rho}(u))\cdot\left(f\otimes\delta\right)=(I_{\nu}(u)f)\otimes\delta.$$
\end{cor}
\begin{proof}
This follows from the previous Lemma. 
\end{proof}
\begin{defn*}
Let $\pr:U(\mathfrak{g}_{\C})\to U(\mathfrak{h}_{\C})$ be
the projection corresponding to the decomposition
$U(\mathfrak{g}_{\C}) = U(\mathfrak{h}_{\C}) \oplus 
\left[
   (\mathfrak{n}_{\C}\oplus\mathfrak{n}_{M})U(\mathfrak{g}_{\C}) +
   U(\mathfrak{g}_{\C})(\bar{\mathfrak{n}}_{\C}\oplus\bar{\mathfrak{n}}_{M}
\right])$
(arising from the decomposition $\mathfrak{g}_{\C}=\mathfrak{n}_{\C}\oplus\mathfrak{n}_{M}\oplus\mathfrak{h}_{\C}\oplus\bar{\mathfrak{n}}_{\C}\oplus\bar{\mathfrak{n}}_{M}$
by the Poincaré-Birkhoff-Witt Theorem).
\end{defn*}
\begin{lem}
\label{lem: z-pr(z) rewrite}For $z\in\mathfrak{Z}^{\leq d}$, we
have\[
z-\pr(z)\in U(\mathfrak{n}_{\C})U(\mathfrak{a}_{\C})^{\leq d-2}U(\mathfrak{k}_{\C}).\]

\end{lem}
\begin{proof}
It suffices to show that $z-\pr(z)\in U(\mathfrak{n}_{\C})U(\mathfrak{g}_{\C})^{\leq d-2}U(\mathfrak{k}_{\C})$,
since $\mathfrak{g}_{\C}=\mathfrak{n}_{\C}\oplus\mathfrak{a}_{\C}\oplus\mathfrak{k}_{\C}$.

Let $\Basis(\mathfrak{n}_{\C})$,
$\Basis(\bar{\mathfrak{n}}_{\C})$,
$\Basis(\mathfrak{n}_{M})$ and
$\Basis(\bar{\mathfrak{n}}_{M})$
be bases for $\mathfrak{n}_{\C}$, $\bar{\mathfrak{n}}_{\C}$,
$\mathfrak{n}_{M}$ and $\bar{\mathfrak{n}}_{M}$, respectively,
consisting of $\mathfrak{h}_{\C}$-eigenvectors. 
Let $\Basis(\mathfrak{a}_{\C})$ and $\Basis(\mathfrak{b}_{\C})$
be bases for $\mathfrak{a}_{\C}$ and $\mathfrak{b}_{\C}$, respectively.

By Poincar{\'e}-Birkhoff-Witt,
one may uniquely express $z$ as a linear combination of terms of the form:
\[
\mathcal{D}=X_{1}\dots X_{n}Y_{1}\dots Y_{m}A_{1}\dots A_{t}B_{1}\dots B_{r}\overline{X}_{1}\dots\overline{X}_{k}\overline{Y}_{1}\dots\overline{Y}_{l}\]
where $X_{*}\in\Basis(\mathfrak{n}_{\C})$, $Y_{*}\in
\Basis(\mathfrak{n}_{M})$,
$A_{*}\in \Basis(\mathfrak{a}_{\C})$, $B_{*}\in\Basis(\mathfrak{b}_{\C})$,$\overline{X}_{*}\in
\Basis(\bar{\mathfrak{n}}_{\C})$
and $\overline{Y}_{*}\in\Basis(\bar{\mathfrak{n}}_{M})$. 
Then $z - \pr(z)$ consists of the sum 
of all terms $\mathcal{D}$ for which $n+m+k+l \neq 0$. 
We show that each
such term satisfies $\mathcal{D}\in
U(\mathfrak{n}_{\C})U(\mathfrak{g}_{\C})^{\leq d-2}U(\mathfrak{k}_{\C})$. 

In view of the fact that $z-\pr(z)$ commutes with $\mathfrak{a}_{\C}$,
one has $n=0$ iff $k=0$. Further, if $n=k=0$,  then the fact that
$z-\pr(z)$ commutes with $\mathfrak{b}_{\C}$ implies $m=0$ iff
$l=0$. Also one has $n+m+t+r+k+l\leq d$. 

We now proceed in a case-by-case
basis, using either the inclusion $\mathfrak{n}_{M}\oplus\mathfrak{b}_{\C}\oplus\bar{\mathfrak{n}}_{M}=\mathfrak{m}_{\C}\subset\mathfrak{k}_{\C}$,
or the observation that for $\overline{X}\in\bar{\mathfrak{n}}_{\C}$ we
have $\theta_{\C}\overline{X}\in\mathfrak{n_{\C}}$, while $\overline{X}+\theta_{\C}\overline{X}\in\mathfrak{k}_{\C}$
(it is $\theta_{\C}$-stable!).
\begin{enumerate}
\item $k=l=0$ is impossible, for this would force $n=m=0$. 

\item $k\geq1$ and $l\geq1$. Then $n\geq 1$ so that $X_{1}\ldots X_{n}\in
U(\mathfrak{n}_{\C})$, $\overline{Y}_{1}\dots\overline{Y}_{l}\in
U(\mathfrak{k}_{\C})$, and $m+t+r+k\leq d-2$.

\item $k=0$ and $l\geq1$. Then $n=0$ and $m\geq1$, so $t\leq d-2$.
Since $\left[\mathfrak{a},\mathfrak{m}\right]=0$ we may commute the
$A$-terms past the $Y$-terms, so that $\mathcal{D}$ is the product
of the $A$-terms (at most $d-2$ of them) and $Y_{1}\dots Y_{m}B_{1}\dots B_{r}\overline{Y}_{1}\dots\overline{Y}_{l}\in U(\mathfrak{k}_{\C})$.
\item $k\geq1$ and $l=0$. Then $n \geq 1$. 
Set $s=Y_{1}\dots Y_{m}A_{1}\dots A_{t}B_{1}\dots B_{r}\overline{X}_{1}\dots\overline{X}_{k-1}$
so that $\mathcal{D}=X_{1}\ldots X_{n}\cdot s\cdot\overline{X}_{k}$.
Since $m+t+r+(k-1)\leq d-1-n\leq d-2$, we have $s\in U(\mathfrak{g}_{\C})^{\leq d-2}$.
Then (recall $\theta_{\C}$ is the complex-linear extension of the
Cartan involution $\theta$ to $\mathfrak{g}_{\C}$), \begin{eqnarray}
\label{eqn: z-pr(z) rewrite}
\mathcal{D}=X_{1}\dots X_{n}s\overline{X}_{k} & = & X_{1}\dots X_{n}\cdot s\cdot(\overline{X}_{k}-\theta_{\C}(\overline{X}_{k}))
\\ & + & X_{1}\dots X_{n}\theta_{\C}(\overline{X}_{k})s 
\nonumber
\\
 & + & X_{1}\dots X_{n}(s\theta_{\C}(\overline{X}_{k})-\theta_{\C}(\overline{X}_{k})s).\nonumber \end{eqnarray}
 From the observation above, the first two terms on the right clearly
belong to $U(\mathfrak{n}_{\C})U(\mathfrak{g}_{\C})^{\leq d-2}U(\mathfrak{k}_{\C})$.
Moreover, $\left[s,\theta_{\C}(\overline{X}_{k})\right]\in U(\mathfrak{g}_{\C})^{\leq d-2}$
(for any $p\in U(\mathfrak{g}_{\C})^{\leq d_{p}},$$q\in U(\mathfrak{g}_{\C})^{\leq d_{q}}$
the general fact $[p,q]\in U(\mathfrak{g}_{\C})^{d_{p}+d_{q}-1}$
follows by induction on the degrees from the formula $[ab,c]=a[b,c]+[a,c]b$).
Thus the third term of (\ref{eqn: z-pr(z) rewrite}) belongs
to $U(\mathfrak{n}_{\C})U(\mathfrak{g}_{\C})^{\leq d-2}U(\mathfrak{k}_{\C})$
also.
\end{enumerate}
\end{proof}
\begin{cor}
\label{cor: z-pr(z) rewrite}Let $z\in \mathfrak{Z}^{\leq d}$.
Then there exists $b=b(z) \in U(\mathfrak{n}_{\C})U(\mathfrak{a}_{\C})^{\leq d-2}$
such that $z-\pr(z)+b(z)\in U(\mathfrak{g}_{\C})\cdot\mathfrak{k}_{\C}$.
\end{cor}

Since $I_{\nu}(\mathfrak{k}_{\C})$ annihilates $\varphi_{0}$ and
$z\cdot\varphi_{0}=\chi_{\nu}(z)\varphi_{0}$,
 we have $I_{\nu}(\chi_{\nu}(z)-\pr(z)+b(z))\cdot\varphi_{0}=0$.
In view of
Corollary \ref{cor: U-on-phi-phi_inf} we obtain:
\begin{equation} 
I_{\nu} \otimes I_{\nu}'(\tau_{\nu+\rho}\pr(z)-\tau_{\nu+\rho}b(z)-
\chi_{\nu}(z))(\varphi_{0}\otimes\delta)=0\label{eq: ann-phi_0-phi_inf}\end{equation}

In what follows, we shall freely identify the algebra $U(\mathfrak{h}_{\C})^{\cW}$
with the Weyl-invariant polynomial functions on 
$\mathfrak{h}_{\C}^{*}$.

Given $\mathcal{P}\in U(\mathfrak{h}_{\C})^{\cW}$,
we denote by 
$\mathcal{P}'\colon\mathfrak{h}_{\C}^{*}\to\mathfrak{h}_{\C}$
its differential. In other words, we identify $\mathcal{P}$
with a polynomial function on $\mathfrak{h}_{\C}^{*}$,
and $\mathcal{P}'$ denotes the derivative of this function;
it takes values in the cotangent space of $\mathfrak{h}_{\C}^{*}$,
which is canonically identified at every point with $\mathfrak{h}_{\C}$.

We shall use the notation $U(\mathfrak{g}_{\C})[\mathfrak{a}_{\C}]^{\leq r}$
to denote polynomials of degree $\leq r$ on $\mathfrak{a}_{\C}^{*}$,
valued in the vector space $U(\mathfrak{g}_{\C})$. 
Note that given $J \in U(\mathfrak{g}_{\C})[\mathfrak{a}_{\C}]^{\leq r}$
 and $\nu \in \mathfrak{a}_{\C}^{*}$ we can 
speak of the ``value of $J$ at $\nu$.'' We denote it by
$J(\nu)$ and it belongs to $U(\mathfrak{g}_{\C})$.

\begin{lem}
\label{lem: calc-J} Let $\mathcal{P}\in U(\mathfrak{h}_{\C})^{\cW}$
have degree $\leq d$. 
Set $H = \frac{\mathcal{P}'(\nu)}{||\nu||^{d-1}} \in \mathfrak{h}_{\C}$. 
Then there exists $J\in U(\mathfrak{g}_{\C})[\mathfrak{a}_{\C}]^{\leq d-2}$
such that \[
I_{\nu} \otimes I_{\nu}'\left(H+\frac{J(\nu)}{\left\Vert \nu\right\Vert
^{d-1}}\right)\cdot\varphi_{0}\otimes\delta=0.\]
(As defined in Section \ref{sec: Notation}, $\left\Vert \nu\right\Vert $ denotes
the norm of $\nu\in\mathfrak{a}_{\C}^{*}$ \wrt the Killing form.)
\end{lem}
\begin{proof}
The map $\gamma_{HC}\colon\mathfrak{Z}\to U(\mathfrak{h}_{\C})^{\cW}$
given by $\gamma_{HC}(z)=\tau_{\rho_{\mathfrak{h}}}\pr(z)$ is an
isomorphism of algebras, the Harish-Chandra homomorphism. With the above identification,
the infinitesimal character of $(V_{K},I_{\nu})$ corresponds to {}``evaluation
at $\nu+\rho-\rho_{\mathfrak{h}}$,'' \ie for
$\mathcal{P}\in U(\mathfrak{h}_{\C})^{\cW}$:
\begin{equation} \label{eq:hchom}
\chi_{\nu}(\gamma_{HC}^{-1}(\mathcal{P}))=\mathcal{\mathcal{P}}(\nu+\rho_{\mathfrak{a}}-\rho_{\mathfrak{h}})
\end{equation}
 (See \cite[Prop 8.22]{Knapp:Rpn_Th_SS_Gps}; \wrt the maximal torus
$\mathfrak{b}_{\C}\subset\mathfrak{m}_{\C}$, the infinitesimal character
of the trivial representation of $\mathfrak{m}_{\C}$ is 
(the Weyl-group orbit of) $\rho-\rho_{\mathfrak{h}}$). 

Given $\mathcal{P}\in U(\mathfrak{h}_{\C})^{\cW}$ of degree $d$,
we set $z=\gamma_{HC}^{-1}(\mathcal{P})$ in (\ref{eq: ann-phi_0-phi_inf}),
writing $b(\mathcal{P})$ for the element $b(z)$. Note that $z\in Z(\mathfrak{g}_{\C})^{\leq d}$,
as the Harish-Chandra homomorphism ``preserves degree'' (see \cite[7.4.5(c)]{Dixmier:EnvlAlg}),
and hence $b(\mathcal{P})\in U(\mathfrak{n}_{\C})U(\mathfrak{a}_{\C})^{\leq d-2}$.

Combining  (\ref{eq: ann-phi_0-phi_inf}) and
(\ref{eq:hchom}):
$\varphi_{0}\otimes \delta$ is then annihilated by the operator
\begin{equation} \label{eq:ann1}
\left( \tau_{\nu+\rho-\rho_{\mathfrak{h}}}\mathcal{P}
-\mathcal{\mathcal{P}}(\nu+\rho_{\mathfrak{a}}-\rho_{\mathfrak{h}})
-\tau_{\nu+\rho}b(\mathcal{P})\right) \varphi_0 \otimes \delta = 0
\end{equation}

Let $\underline{x} = (x_1, \dots, x_n), \underline{y} = (y_1, \dots, y_n)$. 
If a polynomial $p\in\C[\underline{x}${]}
has degree $d$,
$p(\underline{x}+\underline{y})-p(\underline{y})=p'(\underline{y})(\underline{x})+q(\underline{x},\underline{y})$
where $q\in\left(\C[\underline{x}]\right)[\underline{y}]$ has degree
at most $d-2$ in $\underline{y}$, and the derivative $p'(\underline{y})$
is understood to act as a linear functional on $\underline{x}$. 

Applying this to $p = \mathcal{P}, \underline{y} = \nu + \rho 
- \rho_{\mathfrak{h}}$ we see that there exists
$J_1 \in U(\lieg_{\C})[\liea_{\C}]^{\leq d-2}$ with $\mathrm{deg}(J) \leq d-2$ 
and 
\begin{equation} \label{eq:huh}
\tau_{\nu+\rho-\rho_{\mathfrak{h}}}\mathcal{P}-\mathcal{P}(\nu+\rho-\rho_{\mathfrak{h}})
= 
\mathcal{P}'(\nu+\rho-\rho_{\mathfrak{h}})+J_1(\nu)
\end{equation}

Now $b(\mathcal{P}) \in U(\mathfrak{n}_{\C}) \cdot
U(\mathfrak{a}_{\C})^{\leq d-2}$,
so the map $\nu \mapsto \tau_{\nu+\rho}b(\mathcal{P})$
can be regarded as an element
$J_2 \in U(\mathfrak{g}_{\C})[\mathfrak{a}_{\C}]^{\leq d-2}$.
Similarly $
\nu \mapsto \mathcal{P}'(\nu + \rho-\rho_{\lieh}) - \mathcal{P}'(\nu)$
defines an element $J_3 \in U(\mathfrak{g}_{\C})[\liea_{\C}]^{\leq d-2}$. 

Combining these remarks with (\ref{eq:ann1}) and (\ref{eq:huh}), we see
that 
$$(\mathcal{P}'(\nu) + J_1(\nu) + J_2(\nu) + J_3(\nu) )
\varphi_0 \otimes \delta = 0$$

Set $J \equiv J_1 + J_2 + J_3$
and divide by $\left\Vert \nu\right\Vert ^{d-1}$ to conclude. 
\end{proof}
\begin{cor}
\label{cor: P'-inv}Let $\mathcal{P}\in U(\mathfrak{h}_{\C})^{\cW}$.
Notations being as in Definition \ref{def: non-deg-simple}
and Lemma \ref{lem:mucty2},
suppose $\{\psi_n\}$ is conveniently arranged. 
Then $\mu_{\infty}(\varphi_0, \delta)$ is $\mathcal{P}'(\tilde{\nu}_{\infty})$-invariant.
\end{cor}
\begin{proof}
It suffices to verify this for $\mathcal{P}$ homogeneous, say 
of degree $d$. 
Combining
Lemma \ref{lem: calc-J} and Lemma \ref{lem:equiv},
and using the homogeneity of $\mathcal{P}$,
we see that there exists $J \in U(\lieg_{\C})[\liea_{\C}]^{\leq d-2}$ so that
$$\left(\mathcal{P}'(\tilde{\nu}_n) + \frac{J(\nu_n)}
{||\nu_n||^{d-1}}\right)  \cdot \mu_{n}(\varphi_0 \otimes \delta) = 0$$

Here $(\mathcal{P'} + \dots)$ acts on 
$\mu_n(\varphi_0 \otimes \delta)$ according
to the natural action of $U(\lieg_{\C})$ on $\DK'$.
Now fix $g \in \DK$. Let $u \rightarrow u^{t}$
be the unique $\C$-linear anti-involution of $U(\lieg_{\C})$
such that $X^t = -X$ for $X \in \lieg_{\C} \subset U(\lieg_{\C})$. 
Then we have for each $d$ 
\begin{equation} \label{eq:exactde}
\mu_{n}\left(\varphi_{0}\otimes\delta\right)\left(\left(\mathcal{P}'(\tilde{\nu}_{n})-\frac{J^{t}(\nu_{n})}{\left\Vert
\nu_{n}\right\Vert ^{d-1}}\right)g\right)=0.\end{equation}

Note that, as $n$ varies, the quantity
$\left(\mathcal{P}'(\tilde{\nu}_{n})-\frac{J^{t}(\nu_{n})}{\left\Vert
\nu_{n}\right\Vert ^{d-1}}\right)g$ remains in a fixed
finite dimensional subspace of $C^{\infty}(X)_K$.
Further, it converges in that subspace to
$\mathcal{P}'(\tilde{\nu}_{\infty}) g$. 

With these remarks
in mind, we can pass to the limit $n \rightarrow \infty$
in (\ref{eq:exactde}) to obtain
$\mu_{\infty}(\varphi_0 \otimes \delta) (\mathcal{P}'(\nu_{\infty}) g) = 0$, 
\ie $\mathcal{P}'(\nu_{\infty})$ annihilates $\mu_{\infty}$ as required. 
\end{proof}
It remains to show that the subspace 
\begin{equation} \label{eq:sdef}S=\left\{
\mathcal{P}'(\tilde{\nu}_{\infty})\mid\mathcal{P}\in
U(\mathfrak{h}_{\C})^{\cW}\right\} \subset\mathfrak{h}_{\C}\end{equation}
contains $\mathfrak{a}_{\C}$. By the Corollary this will show that
$\mathfrak{a}$ annihilates any limit measure, or that this measure
is $A$-invariant.

\begin{lem} \label{lem:previous}
Let $W_{0}\subset\cW$ be the stabilizer of $\tilde{\nu}_{\infty}
 \in\mathfrak{a}_{\C}^{*}$, and define $S$ as in (\ref{eq:sdef}). 
 Then $S=\mathfrak{h}_{\C}^{W_{0}}$.
In particular, if $\tilde{\nu}_{\infty}$ is regular,
then $S$ contains $\liea_{\C}$. 
\end{lem}
\begin{proof}
This can be seen either from the fact that $S$ is the image of the
map on cotangent spaces induced by the quotient map $\mathfrak{h}_{\C}^{*}\rightarrow\mathfrak{h}_{\C}^{*}/W_{0}$,
or more explicitly: first construct many elements in $U(\mathfrak{h}_{\C})^{\cW}$
by averaging over $\cW$, and then directly compute derivatives to
obtain the claimed equality.

$W_0$ is generated
by the reflections in $\cW$ fixing $\tilde{\nu}_{\infty}$. In the case where
$\tilde{\nu}_{\infty}$ is regular as an element in $i\mathfrak{a}_{\R}^{*}$,
the corresponding roots must be trivial on all of $\mathfrak{a}_{\C}^{*}$.
In particular, any element of $W_{0}$ fixes all of $\mathfrak{a}_{\C}$.
\end{proof}

\begin{cor}
\label{cor: A-inv} 
Let notations be as in Proposition \ref{prop: equivariance}. 
Then any weak-{*} limit $\sigma_{\infty}$ of the measures $\sigma_n$ is $A$-invariant. 
\end{cor}
\begin{proof}
After passing to an appropriate subsequence,
we may assume that $\{\psi_n\}$ is conveniently arranged. 
Proposition \ref{prop: equivariance}, (\ref{claim1}),  shows that
$\sigma_{\infty}(g)  = \mu_{\infty}(\varphi_0, \delta) (g)$ whenever
$g \in \DK$. 
Corollary \ref{cor: P'-inv} and Lemma \ref{lem:previous}, together
with the fact that $\DK$ is dense in $C_0(X)$,
show that $\sigma_{\infty}$ is $A$-invariant. 
\end{proof}

\section{Complements}
In this section we gather together several points complementing
the main text. 

\subsection{\label{sub: extend}Extensions to general $G$}

In practice we wish to apply our result to groups which are
slightly more general that the ones considered above.  
Here we briefly discuss extensions of the present work to reductive groups. 

From now on let $G$ be a linear connected reductive Lie group,  $\Gamma<G$ a
lattice.  (For ``linear connected reductive'', we 
follow the definition of \cite[Chapter 1]{Knapp:Rpn_Th_SS_Gps}.)
We set $X=\Gamma\backslash G$ as before, and define in addition
$X_Z=Z\Gamma\backslash G$, with $Z=Z(G)$.  $X_Z$ has finite volume
\wrt the $G$-invariant measure.  In a similar fashion we shall consider
$Y=X/K$ and $Y_Z = X_Z/K$.

Since $G$ is linear connected reductive, we have a decomposition
$\mathfrak{g}=\mathfrak{z}\oplus_{j}\mathfrak{g}^{(j)}$ where
$\mathfrak{z}=Z_{\mathfrak{g}}$, and each $\mathfrak{g}^{(j)}$ is a simple Lie
algebra (an orthogonal decomposition \wrt the Killing form), leading to a
decomposition $G=Z_{G}\times\prod_{j}G^{(j)}$ (almost direct product),
where the $G^{(j)}$ are connected semi-simple or compact normal
subgroups.

Choosing the Cartan involution, the subalgebra $\liea$, etc.\ compatible
with this decomposition, let $K=K_Z\times\prod_{j}K^{(j)}$
be the $\Theta$-fixed maximal compact subgroup. If $G^{(j)}$ is
compact then $K^{(j)}=G^{(j)}$, of course.
Note that the subgroup $M=Z_{K}(\mathfrak{a})$ now includes
the compact part of the center, as well as all compact factors.

For a unitary character $\omega\in\hat{Z}$, let $L^{2}(X,\omega)$
denote the space of all measurable $f:X\to\C$ such that
$f(zg)=\omega(z)f(g)$
for all $z\in Z$, and such that
$\left\Vert f\right\Vert^2 = \int_{X_{Z}}|f(x)|^2 dx<\infty$.
If $\omega$ is \emph{unramified} (\ie trivial on $Z(G)\cap K$),
then set $L^{2}(Y,\omega) = L^2(X, \omega)^K$.
If $\omega$ is unramified
and $\psi \in L^2(Y,\omega)$, then $|\psi(y)|^2$ is $Z$-invariant,
and we can define a finite measure $\bar{\mu}_\psi$ on $Y_Z$ as before.

An eigenfunction $\psi\in L^{2}(Y,\omega)$ still generates an
irreducible subrepresentation of $G$ in $L^2(X,\omega)$.
From this we obtain, as in Section \ref{sec: Lift},
a norm-reducing intertwining operator
$R_{\psi} \colon(V_{K},I_{\nu}) \to L^{2}(X,\omega)$, and 
(as in Definition \ref{def:mudef}) a map
$\mu_{\psi}\colon V_K\otimes V_K'\to \left(\Cic(X_Z)_K\right)'$ as before
(note that for $f_{1},f_{2}\in V_{K}$, $R_{\psi}(f_{1})\overline{R_{\psi}(f_{2})}$
is $Z$-invariant since $\omega$ is unitary, and as before its $L^1$ norm
is at most the product of the $L^2$ norms of $f_1,f_2\in V$).

Let $\{\omega_n\}_{n=1}^{\infty}$ be a sequence
of unramified characters of $Z$.
We now consider a sequence of eigenfunctions
$\left\{ \psi_{n}\right\}_{n=1}^{\infty}$ such that
$\psi_n\in L^{2}(Y,\omega_n)$, with intertwining operators $R_{n}$
and parameters $\nu_{n}\in\mathfrak{a}_{\C}^{*}$, and assume that the
$\nu_{n}$ escape to infinity.

\begin{defn}
\label{def: non-deg-red}Call the sequence \emph{non-degenerate}
if for every non-compact $j$, the sequence
$\{\nu_n^{(j)}\}\subset \left(\liea_\C^{(j)}\right)^*$ is
non-degenerate in the sense of Definition \ref{def: non-deg-simple}.
\end{defn}
\begin{rem}
As before, for a non-degenerate sequence we have $\mathrm{Re}(\nu_{n})=0$
and $\textrm{Im}(\nu_{n})$ regular for large enough $n$. However, the
rates at which the different components of $\nu_n$ tend to infinity 
need not be the same. 

Indeed, defining a $\tilde{\nu}_{n}^{(j)}$ for each $j$ by normalizing $\nu_n^{(j)}$,
and passing to a subsequence where they all converge,
the non-degeneracy assumption amounts to assuming that the
limits $\tilde{\nu}^{(j)}$ are regular (\ie do not lie on any wall).
Of course, the rate of convergence at different $j$ may be different.
\end{rem}
Lemma \ref{lem: int-parts} and its subsequent Corollary continue to hold
(replace $\D$ with $\Cic(X_Z)$).  The only modification to the proof is
that one should only consider functions $p_X$ given by
$X\in\mathfrak{g}^{(j)}$, rescaling by $\left\Vert \nu_{n}^{(j)}\right\Vert$.
The Stone-Weierstrass argument will show that the algebra generated by
these ``limited'' $p_{X}$ is dense.
Defining the lift as before (using the $\delta$ distribution at $1\in\mk$),
we obtain the positivity of the limits.

In the same vein it is clear that by using $U(\mathfrak{g}_{\C}^{(j)})$
and its center (which is contained in the center of $U(\mathfrak{g}_{\C})$),
the analysis of Section \ref{sub: A-inv} shows that a non-degenerate
limit is $\mathfrak{a}^{(j)}$-invariant for all non-compact $j$,
and hence $A$-invariant. As before, every $\mu_{R_{n}}$ is $M$-invariant,
hence so is $\mu_{\infty}$.

\subsection{Degenerate limits.}
It is an interesting and natural problem to extend the results
of the present paper to degenerate limits, \ie
sequence of eigenfunctions $\psi_n$ such
that $\frac{\nu_n}{\Vert\nu_n\Vert}$ converges
to one of the walls of a Weyl chamber. 

The non-degeneracy assumption was used in several places in the above
arguments.  The first was in the assertion that the intertwining maps
from the models $(\pi_\nu, V_K)$ to $L^2(X)$ were isometries for the $L^2$ norm
on $V_K$, so that the total variation of the measures $\mu^T_\psi(f_1,f_2)$
was bounded independently of the parameter $\nu$ of $\psi$.
Secondly, we used it in the proof of positivity of the limit measures by
integration by parts. Finally, it was used to conclude
that the limit measures is indeed invariant under the full
Cartan subgroup $A$.

The first use can be removed in a straightforward manner resulting in a
lift of the limit measure which is a positive measure on $X$. However,
the question of invariance is more subtle, and one might expect the methods
presented here to only show invariance under an appropriate subtorus of $A$.
We hope to revisit this issue in the future.

\subsection{Geometry of the Cartan flow and flats.} \label{chamber}

A symmetric space comes with a rich structure of flat subspaces; these
are an important part of the large-scale geometry of the space. Our
aim here is to discuss the connection of the Cartan flow (i.e.
the action of $A$ on $X/M$)
with the structure of flats. Crudely speaking, the Cartan flow is
analogous to the geodesic flow, but with {}``geodesic'' replaced
by {}``flat.''  This highlights the fact that the present result
is a generalization of the rank $1$ situation, where flats
\emph{are} geodesics. (The present result, however, is new
even in the case of hyperbolic $3$-space,
on account of its equivariance.)

Let $G$ and other notations be as fixed in Section \ref{sub: notation},
and let $r$ be the real rank of $G$, $W=\rW$ the Weyl group.

An \emph{$r$-flat} in $S$ is, by definition, a subspace isometric
to $\R^{r}$ with the flat metric. Given any $r$-flat $F\subset S$
and a point $P\in F$, there is a canonical $W$-conjugacy class $\mathcal{C}_{P,F}$
of isometries from $\mathfrak{a}$ to $F$, all mapping $0$ to $P$.
Indeed, we may assume that $P=x_{K}$, in which case we may identify
$F$ (via the inverse exponential mapping) with a subset of $\mathfrak{p}$,
which may be shown (see \cite{Mostow:Rigidity}) to be a maximal abelian
subspace. In particular, this subset is conjugate under $K$ to $\mathfrak{a}$,
and this conjugacy is unique up to the action of the Weyl group, whence
the assertion.

An \emph{orientation} $\varphi$ of the pair $(P,F)$ will be an element
$\varphi\in\mathcal{C}_{P,F}$; there are therefore precisely $|W|$
orientations for any pair $(P,F)$. A \emph{chamber} will be a triple
$(P,F,\varphi)$ of a point $P$, a flat $F$ containing $P$, and
an orientation for $(P,F,\varphi)$. In the case $r=1$, a chamber
is equivalent to a geodesic ray: given a chamber $(P,F,\varphi)$,
the set $\varphi^{-1}([0,\infty))$ is a geodesic ray beginning at
$P$.

The \emph{chamber bundle} of $S$, denoted $\mathcal{C}S$, will be
the set of all chambers. $G$ acts transitively on $\mathcal{C}S$
and the stabilizer of a point is conjugate to $Z_{K}(\mathfrak{a})$ ($=M$).
In particular, $\mathcal{C}S$ has the structure of a differentiable
manifold, and it is a fiber bundle over $S$; each fiber is isomorphic
to $K/Z_{K}(\mathfrak{a})$.

The additive group of $\mathfrak{a}$ acts in an evident way on $\mathcal{C}S$:
given $X\in\mathfrak{a}$ and a chamber $(P,F,\varphi)$, one defines
$X(P,F,\varphi)=(\varphi^{-1}(X),F,\varphi')$, where there is a unique
choice of $\varphi'$ that makes this a continuous action. In particular,
$\mathcal{C}S$ carries a natural $\R^{r}$ action. In the case $r=1$
this is the geodesic flow on the unit tangent bundle.

Finally, if $\Gamma$ is any discrete subgroup of $G$, one sees that
$\R^{r}$ acts on $\Gamma\backslash\mathcal{C}S$, which fibers over
$\Gamma\backslash S$. The main result of the present paper may be
phrased as follows: a measure on $\Gamma\backslash S$ arising from
a limit of eigenfunction measures lifts to an $\R^{r}$-invariant
measure on $\Gamma\backslash\mathcal{C}S$.

\subsection{Relation to $\Psi$DOs.}
Zelditch's original proof for hyperbolic $2$-space
involved the construction of an equivariant pseudodifferential
calculus based on the non-Euclidean Fourier transform of Helgason. 
It is certainly reasonable to expect that this could be generalized
to higher rank; however, for the application to quantum chaos,
the methods of this paper seem more efficient. 
In either approach, the positivity and Cartan invariance
require proof. 

Of course, the two methods are very closely linked. 
In this section we translate the representation-theoretic methods
of this paper to the microlocal viewpoint. 
In fact, we will only do the bare minimum to show
that the microlocal lifts constructed in the present
paper are ``compatible'' with the standard construction for a general
Riemannian manifold described in \cite{CdV:Avg_QUE}.
We will also only sketch the proof; it is more or less formal. 

From the microlocal viewpoint the system under consideration resembles
completely integrable systems, in that there are several commuting observable;
see \eg \cite{TothZelditch:CompletelyInt}.
Of course, the Cartan flow differs from the completely integrable case
in that it is very chaotic. 

Initially let the notation be as in the introduction; in particular
let $\Man$ be a compact Riemannian manifold, $\Delta$ the Laplacian
on $\Man$, $S^{*}\Man$ the unit cotangent bundle. 
We fix a quantization
scheme $\Op$ that associates to a smooth function $a$ on $S^{*}\Man$
a pseudo-differential operator $\Op(a)$ on $\Man$ of order $0$. Let
$\psi_{n}$ be a sequence of eigenfunctions of $\Delta$ with eigenvalues
$\lambda_{n}\rightarrow-\infty$, s.t. the measure $\bar{\mu}_{\infty}=\lim_{n\to\infty}|\psi_{n}|^{2}d\rho$
exists. Then, after possibly passing to a subsequence, the limit $a\mapsto\langle\Op(a)\psi_{n},\psi_{n}\rangle$
exists for all $0$-homogeneous $a$ and defines a positive measure $\mu_{\infty}$
that lifts $\bar{\mu}_{\infty}$. We shall refer to this
as a \emph{standard microlocal lift}.

Now let us follow the notation of Section \ref{sub: notation}. For
simplicity we shall assume $G$ simple and center-free and $\Gamma\subgp G$
co-compact. We shall also identify $\mathfrak{g}$ and $\mathfrak{p}$
with their duals by means of the Killing form, and we will identify
the tangent and cotangent bundle of $Y$ by means of the Riemannian
structure (induced from the Killing form as well). We denote by $\left\Vert \cdot\right\Vert $
the norm induced on $\mathfrak{g}$ and $\mathfrak{g}^{*}$ by the
Killing form.

Let us recall more carefully the connection between $X$ and the tangent
bundle of $Y$. As before set $X=\Gamma\backslash G$, $Y=\Gamma\backslash G/K$,
$S=G/K$, and let $\pi\colon X\to Y$ denote the natural projection.
Let $TS$ and $TY$ denote the tangent bundles of $S$ and $Y$, and
let $x_{K}\in S$ be the point with stabilizer $K$.  Let $T^{1}Y \subset TY$
be the unit tangent bundle; we will often implicitly
identify functions on $T^{1}Y$ with $0$-homogeneous functions on $TY$,
and in particular functions on $T^{1}Y$ gives rise to pseudodifferential
operators of order $0$. 

We shall endow
$G\times\mathfrak{p}$ with the left $G$-action given by $g(h,Y)=(gh,Y)$,
and with the right $K$-action given by $(h,Y)k=(hk,k^{-1}Yk)$. There
is a natural map $G\rightarrow S$ given by $g\mapsto gx_{K}$. This
lifts to a $G$-equivariant map $G\times\mathfrak{p}\rightarrow TS$;
this latter map is specified by requiring that its restriction to
$\{ e\}\times\mathfrak{p}$ be the usual identification of $\mathfrak{p}$
with the tangent space to $S$ at $x_{K}$. Taking quotients by $\Gamma$,
we descend to a map also denoted $\pi\colon X\times\mathfrak{p}\rightarrow TY$.
This map is constant on $K$-orbits, and factors through to a map
$X \times \mathfrak{p}/K \rightarrow TY$. 

In view of our identification of tangent and cotangent bundles, the
symbol of a pseudodifferential operator on $Y$ may then be regarded
as a $K$-invariant function on $X\times\mathfrak{p}$.
We shall fix
a quantization scheme $\Op$ that associates to such a symbol a pseudo-differential
operator on $Y$.

Let $\left\{ \psi_{n}\right\} _{n=1}^{\infty}\subset L^{2}(Y)$ be
a sequence of eigenfunctions on $Y$ with parameters $\nu_{n}\in \mathfrak{a}^{*}$
and so that $\frac{\nu_{n}}{\left\Vert \nu_{n}\right\Vert }\rightarrow\tilde{\nu}$.
We shall assume that $\left\{ \psi_n \right\}$ is conveniently arranged
in the sense of Definition \ref{def: non-deg-simple}. 
We let $\lambda_{n}$ be the Laplacian eigenvalue of $\psi_{n}$ (this differs
by a constant from $-\left\Vert \nu_{n}\right\Vert ^{2}$, in fact). We can and will
also regard the $\psi_{n}$ as $K$-invariant functions on $X$. Associated
to each $\psi_{n}$ is an $G$-intertwiner $R_{n}:(V_{K},I_{\nu_{n}})\rightarrow L^{2}(X)$.

We shall use $o(1)$ to denote quantities with
go to $0$ as $\left\Vert \nu_{n}\right\Vert \to\infty$.

Let other notations be as in Section \ref{sub: positivity}. 
The relation between the $\Psi D O$ viewpoint and the methods of
this paper are summarized in:
\begin{prop*}
Let $a\in C^{\infty}(T Y)$ be such that $a$ is $0$-homogeneous.
Let $g \in C^{\infty}(X)$ be defined by $g(x) = a(\pi(x,\tilde{\nu}))$. 
Suppose that $g$ is right $K$-finite. 
Then\begin{equation}
\left\langle \Op(a)\psi_{n},\psi_{n}\right\rangle =\mu_{n}(\varphi_0,
\delta)(g)+o(1).\label{eq: cdv}\end{equation}

\end{prop*}
It follows that if $\mu_{\infty,T^{1}Y}$ is a standard microlocal
lift then $\mu_{\infty,T^{1}Y}$ is supported on 
$\pi(X\times\{\tilde{\nu}\})$, and the restriction
of $\mu_{\infty,T^{1}Y}$ to this copy of $X$ is a microlocal lift
in the sense of the current paper.

\begin{proof}
In three stages.

\emph{First step.}
We first verify that, if $g=0$, then
$\langle\Op(a)\psi_{n},\psi_{n}\rangle=o(1)$.

Let $P$ be a $K$-invariant polynomial on $\mathfrak{p}$ of degree
$d$ and consider the function $\tilde{P}:(x,A)\in X\times\mathfrak{p}\mapsto P(A)$.
The function $\tilde{P}$ descends to $TY$, 
and there is an invariant differential operator $\mathcal{D}_{P}$ 
on $Y$ of degree $d$ whose symbol agrees with $\tilde{P}$. 
Since $\psi_{n}$ is an eigenfunction for the ring of invariant differential
operators, it follows in particular that
$\psi_{n}$ is an eigenfunction for $\mathcal{D}_{P}$ with eigenvalue
$P(\nu_{n})$. It follows that, for any $b\in C^{\infty}(X \times
\mathfrak{p})^K$, \begin{equation} \label{eq:doon}
\frac{P(\nu_{n})}{\left\Vert \nu_{n}\right\Vert ^{d}}\langle\Op(b)\psi_{n},
\psi_{n}\rangle=\frac{\langle\Op(b)\mathcal{D}_{P}\psi_{n},\psi_{n}\rangle}
{\left\Vert \nu_{n}\right\Vert ^{d}}
=\frac{\langle\Op(b\tilde{P})\psi_{n},\psi_{n}\rangle}
{\left\Vert \nu_{n}\right\Vert ^{d}}+o(1).\end{equation}

(\ref{eq:doon}) implies, in particular, that if $P(\tilde{\nu}) = 0$
the statement of the Proposition holds for $a = b \tilde{P}$.
We can deduce the claim of the first step by density:
if $a|_{X \times \{\tilde{\nu}\}^K}$ is identically $0$,
then one can verify that $a$ may be densely approximated (in the topology
induced by symbol-norm) by linear combinations of functions
$b \cdot P$ where $P(\tilde{\nu}) = 0$. 
We conclude using $L^2$-bounds on pseudodifferential operators 
(\cite[Thm.\ 18.1.11]{Hormander:PD_Ops_III} and remarks after proof.)

\emph{Second step.} We next construct an explicit class of test functions $a$ for which
(\ref{eq: cdv}) holds.

Let $\sigma\in C^{\infty}(X)_{K}$, $u\in U(\mathfrak{g})$ 
of degree $\leq d$, and let $\pi^{*}$ and $\pi_{*}$ be, respectively, the pull-back
and push-forward operations on functions arising from $\pi:X\to Y$.
(In other words, $\pi_{*}$ is obtained by integrating along $K$-orbits.) 
Let $\mathrm{mult}_{\sigma}$ be the operation {}``multiplication
by $\sigma$'' on $C^{\infty}(X)$. We can define by the spectral
calculus of self-adjoint operators an endomorphism $(1-\Delta)^{-d/2}\colon C^{\infty}(Y)\rightarrow C^{\infty}(Y)$.
We then define an endomorphism of $C^{\infty}(Y)$ via the rule \[
\mathrm{MyOp}(\sigma):f\mapsto\pi_{*}\circ\mathrm{mult}_{\sigma}\circ
u\circ\pi^{*}\circ(1-\Delta)^{-d/2}f\]

In other words, one applies $(1-\Delta)^{-d/2}$, lifts the resulting
function to $X$, applies $u$ and multiplies by $\sigma$, and pushes
back down to $Y$.

Regard $u$ as defining (its ``symbol'') a polynomial function $u_{d}$ of degree $d$
on $\mathfrak{g}^{*}$ (therefore on $\mathfrak{g}$) and let $a_{\sigma,u}$
be the following $K$-invariant function on $X\times\mathfrak{p}$:
\[
a_{\sigma,u}:(x,A)\in X\times\mathfrak{p}\mapsto\frac{1}{||A||^{d}}\int_{K}\sigma(xk)u_{d}(k^{-1}Ak)dk\]

We verify (\ref{eq: cdv}) for $a=a_{\sigma,u}$ and $g_{\sigma,u}\stackrel{\mathrm{def}}{=}a_{\sigma,u}|X\times\{\tilde{\nu}\}$.
Note that \begin{equation} 
g_{\sigma,u}(x)=\int_{K}\sigma(xk)u_{d}(k^{-1}\tilde{\nu}k)dk\label{geqn}
\end{equation}

The operator $\mathrm{mult}_{\sigma} \circ u$ is
clearly a differential operator on $X$,
and one deduces that the operator $\pi^{*}  \circ \mathrm{mult}_{\sigma}
\circ u \circ \pi^{*}$ is, in fact, a differential operator
on $Y$.  One computes that the symbol
of this latter operator is associated
to the $K$-invariant function $(x,A)\in
X\times\mathfrak{p}\mapsto\int_{K}\sigma(xk)u_{d}(k^{-1}Ak)dk$.
We deduce that: \begin{equation}
\langle\mathrm{MyOp}(\sigma)\psi_{n},\psi_{n}\rangle-\langle\mathrm{Op}(a_{\sigma,u})\psi_{n},\psi_{n}\rangle=o(1)\label{eqnone}\end{equation}
 Further, if we regard $\psi_{n}$ as a $K$-invariant function on
$X$: \begin{equation}
\langle\mathrm{MyOp}(\sigma)\psi_{n},\psi_{n}\rangle=(1-\lambda_{n})^{-d/2}\int_{X}\overline{\psi_{n}(x)}\sigma(x)(u\psi_{n})(x)\, dx\label{eqntwo}\end{equation}

On the other hand, recall the definition of $\mu_n$
from Section \ref{sub: def-lift}.
Let $\delta_{(N)}$ be the $N$-truncation of $\delta$ (see
Definition \ref{def: vdef}). Choosing $N$ sufficiently large, 
we have $\mu_n(\varphi_0, \delta) g_{\sigma, u} 
= \mu_n^T (\varphi_0, \delta_N) g_{\sigma,u}$; in particular
\begin{eqnarray}
\mu_{n}(\varphi_0,\delta)g_{\sigma,u} & = & 
\int_{X}dx\,\psi_{n}(x)\cdot \overline{R_{n}(\delta_{N})(x)}\cdot\int_{K}\sigma(xk)u_{d}(k^{-1}\tilde{\nu}k)dk\\
 & = & \int_{X}\int_{K}dx \, dk\,
 \psi_{n}(x) \overline{R_{n}(\delta_{N})(xk^{-1})}\sigma(x)u_{d}(k^{-1}\tilde{\nu}k)\\
 & = & \int_{X}dx \, \psi_{n}(x)\sigma(x)\overline{R_{n}\left(\int_{K}dk\,
 \overline{u_{d}(k\tilde{\nu}k^{-1})}I_{\nu_{n}}(k)\,\delta_{N}\right)}\label{eqnthreeminusone}\end{eqnarray}

At the last step, we make the substitution
$k \mapsto k^{-1}$, 
and use the fact that the representation
$I_{\nu_{n}}|_{K}$ is just the operation
of right translation.

To simplify this further, we use Lemma \ref{lem: mathfrak{g}-on-ind}.

Let $p_{u}$ be the function $k\mapsto u_{d}(k\tilde{\nu}k^{-1})$;
it defines a function on $M\backslash K$ and thus we can regard 
$p_{u}\in V_{K}$.  
Denote by $\overline{p_u}$ the complex conjugate of $p_u$.
Since $\delta_{N}$ is, as a function on $M\backslash K$, an approximation
to a $\delta$-function, we have as $N \rightarrow \infty$:
$\int_{K} dk\, \overline{u_{d}(k\tilde{\nu}k^{-1})}I_{\nu_{n}}(k)\delta_{N}
\rightarrow p_{u}$. Here 
the convergence occurs in $C(M\backslash K)$.
It follows that for any $n$:
\begin{equation}
\label{eqn:hungry}\mu_n(\varphi_0, \delta) g_{\sigma, u}
= \langle \psi_n(x) \sigma(x), R_n(\overline{p_u}) 
\rangle_{L^2(X)}\end{equation}

In view of the definitions, the right-hand side 
of (\ref{eqn:hungry}) is just $\mu_n(\varphi_0, \overline{p_u}) (\sigma)$. 
By (the proof of) Lemma \ref{lem: int-parts}, 
$\mu_n(\varphi_0, \overline{p_u})(\sigma) = \mu_n(p_u \varphi_0,
\varphi_0)(\sigma) + o(1)$. Consequently,
\begin{equation} \label{eqn:hungry2}
\mu_n(\varphi_0, \delta) g_{\sigma,u} = \int_{X} 
\overline{\psi_n(x)} R_n(p_u)(x) \sigma(x) dx + o(1).
\end{equation}
On the other hand, a computation with 
Lemma \ref{lem: mathfrak{g}-on-ind}  shows that
$$p_u - \frac{ I_{\nu_n}(u)  \varphi_0}{(1-\lambda_n)^{d/2}}
\rightarrow 0 \mbox{ in } L^2(M \backslash K).$$
Combining this with (\ref{eqn:hungry2}), we obtain:
\begin{equation} \label{eqnthree}\mu_n(\varphi_0, \delta) g_{\sigma,u} - 
(1-\lambda_n)^{-d/2} \int_X \overline{\psi_n(x)} \sigma(x) u \psi_n(x) dx 
= o(1).\end{equation}

In view of (\ref{eqnone}), (\ref{eqntwo}) and (\ref{eqnthree})
we have verified (\ref{eq: cdv}) in the case of $a=a_{\sigma,u}$.

\emph{Third step.} 
Note that, in the statement of the Proposition,
the function $g$ is necessarily right $M$-invariant. 
In view of what has been proved, it now suffices to check that functions
of the form $g_{\sigma,u}$ (see (\ref{geqn}))
span $C^{\infty}(X/M)_{K}$. This is easily
reduced to checking that the linear span
of the functions $k\mapsto u_{d}(k^{-1}\tilde{\nu}k)$
is a dense subspace of $C(M\backslash K)$.
This is shown in the proof of Lemma \ref{lem: int-parts}.
\end{proof}

\bibliographystyle{amsplain}
\bibliography{SQUE_Lift}

\end{document}